\documentclass[a4paper,fleqn]{cas-sc}
\usepackage{amsmath,amssymb,amsthm,mathrsfs}
\usepackage{graphicx,xcolor}
\DeclareRobustCommand{\HW}[1]{{#1}}
\ExplSyntaxOn
\RenewDocumentCommand \printorcid { } { }
\cs_set:Npn \__first_footerline: { \null }
\ExplSyntaxOff
\usepackage[font=footnotesize]{subfig}
\usepackage{placeins}
\usepackage[numbers,sort&compress]{natbib}
\hypersetup{colorlinks=true,linkcolor=blue,citecolor=blue,urlcolor=blue}
\newtheorem{theorem}{Theorem}[section]
\newtheorem{lemma}[theorem]{Lemma}
\newtheorem{proposition}[theorem]{Proposition}
\theoremstyle{remark}
\newtheorem{remark}[theorem]{Remark}
\numberwithin{equation}{section}

\begin{document}
\shorttitle{Learning the Boltzmann Collision Operator}
\shortauthors{H. Wang and Z. Cai}
\title[mode=title]{Learning the Boltzmann Collision Operator with Structural Constraints for Fast Kinetic Simulation}
\author[1]{Hao Wang}
\ead{haowang@u.nus.edu}
\credit{Conceptualization, Methodology, Software, Formal analysis, Visualization, Writing -- original draft, Writing -- review \& editing}
\author[1]{Zhenning Cai}
\ead{matcz@nus.edu.sg}
\credit{Conceptualization, Methodology, Supervision, Writing -- review \& editing}
\affiliation[1]{organization={Department of Mathematics, National University of Singapore},postcode={119076},country={Singapore}}

\begin{abstract}
The high-dimensional integration in the Boltzmann collision operator is a major computational bottleneck in deterministic kinetic simulations. We develop a structure-preserving neural operator that provides an efficient approximation of this nonlinear collision mapping. The operator maintains the gain--loss structure of the original Boltzmann collision operator, and, by construction, is proportional to the square of the gas density. Conservative moment projection and entropy correction enforce conservation of mass, momentum, and energy, as well as entropy dissipation. By construction, the collision term also vanishes at all Maxwellians. In the numerical scheme for simulations of gas dynamics, we establish nonnegativity of the distribution function with the learned collision operator and prove exponential relaxation to the moment-matched equilibrium when a positive spectral gap is imposed. Numerical experiments for a constant collision kernel demonstrate a favorable accuracy--cost tradeoff relative to Fourier spectral evaluations. %
Although trained on smooth distributions, the same operator handles discontinuous distribution functions and is embedded without retraining in steady Couette flow, normal-shock, and two-dimensional lid-driven cavity solvers. These results demonstrate a reusable approximation for repeated collision evaluations that retains kinetic structure and supports accurate deterministic simulation.
\end{abstract}

\maketitle

\section{Introduction}
\label{sec:introduction}

Kinetic models play a fundamental role in describing dilute gases and nonequilibrium transport phenomena, with applications ranging from rarefied gas dynamics and aerospace engineering to microscale and nanoscale flows and vacuum systems~\cite{cercignani2002boltzmann,shen2005rarefied}. When the molecular mean free path is not negligible compared with the characteristic length scale, continuum descriptions based on the Navier--Stokes equations may no longer adequately capture the underlying dynamics~\cite{torrilhon2016modeling}. The Boltzmann equation provides a fundamental kinetic description through the evolution of the particle distribution function $f$ in phase space~\cite{Cercignani1988}. Its numerical solution remains challenging, primarily because deterministic solvers require repeated evaluation of the nonlinear collision operator over the spatial domain.

The main computational challenge arises from the high-dimensional integral structure of the Boltzmann collision operator~\cite{DimarcoPareschi2014}. Fourier spectral methods provide an accurate deterministic framework for its approximation, but direct evaluation in three-dimensional velocity space requires $\mathcal{O}(N^6)$ operations~\cite{original_fourier_collision}, where $N$ denotes the number of Fourier modes in each velocity direction. Fast spectral algorithms exploit the Fourier-space convolution structure together with fast Fourier transforms to reduce the computational complexity to \(\mathcal{O}(MN^4\log N)\) for general collision kernels~\cite{fast_spectral_method_boltzmann1,fast_spectral_method_boltzmann2}, where \(M\) denotes the number of quadrature points for angular integration over \(\mathbb{S}^2\). Even at this reduced complexity, repeated collision evaluations remain expensive in spatially inhomogeneous simulations, where the operator must be computed at every spatial location.

This cost has motivated a range of alternative strategies. The direct simulation Monte Carlo (DSMC) method replaces deterministic collision integration with stochastic particle collisions and has become a standard tool for rarefied gas simulations~\cite{DSMC}, although statistical noise can limit its efficiency when high accuracy is required~\cite{DSMC_error}. Simplified relaxation models, such as the BGK and Shakhov models~\cite{BGK1,BGK2}, offer much lower computational cost but provide only an approximate representation of the full collision dynamics. For steady-state problems, accelerated iterative methods use the linearized collision operator to construct efficient update directions, whose Fourier spectral evaluation for variable hard-sphere gases can be performed with \(\mathcal{O}(N^4\log N)\) complexity~\cite{linearized_boltzmann_collision_operator,steady_boltzmann1}. These strategies reduce the computational burden through stochastic sampling, model simplification, or accelerated convergence.

Learning the full nonlinear collision mapping offers a complementary way to reduce the cost of repeated evaluations. Neural operators provide a flexible framework for learning nonlinear mappings between function spaces directly from data~\cite{FNO,Kovachki2023NeuralOperator,deeponet}. Hybrid scientific machine learning retains reliable numerical components while learning only computationally expensive or unresolved parts of the dynamics, as recently demonstrated by alternating neural integrators that couple learned components with existing simulators without modifying the underlying solver~\cite{Willard2023,Wang2026Learning}. The Boltzmann equation is particularly well suited to this paradigm. Standard numerical methods can continue to handle transport, while the expensive mapping from $f$ to $\mathcal{Q}[f,f]$ can be replaced by a reusable learned surrogate.

Within kinetic modeling, several learned approaches have targeted the Boltzmann collision operator directly. RelaxNet uses a parameterized relaxation form to preserve positivity, while both conservation and the $H$-theorem hold only up to training error~\cite{RelaxNet}. FourierSpecNet combines neural learning with a Fourier spectral representation for efficient and resolution-independent approximation~\cite{FourierSpecNet}. Approximation accuracy alone is not sufficient for reliable collision dynamics. They must also respect positivity, conservation laws, entropy dissipation, and Maxwellian equilibria. This issue is not specific to learned surrogates. In deterministic Boltzmann solvers, these requirements have motivated schemes that preserve entropy and positivity~\cite{CaiFanYing2018}, spectral methods that preserve collision moments~\cite{PareschiRey2022}, and spectral approximations that preserve both positivity and moments~\cite{CaiLinLin2024}. Highly accurate discretizations can still violate positivity or the $H$-theorem~\cite{CaiFanYing2018}, while loss of collision invariants may lead to incorrect long-time behavior~\cite{PareschiRey2022}. This motivates learned collision operators that enforce kinetic structure directly rather than relying solely on approximation accuracy or training-based regularization.

In this work, we develop a neural operator approximation of the full Boltzmann collision operator with exact preservation of key kinetic structure. For admissible distributions, the construction preserves mass, momentum, and energy, dissipates entropy, and is consistent with Maxwellian equilibria independently of the trained parameters. We also establish positivity of the continuous collision dynamics under suitable regularity and integrability assumptions and derive a sufficient time step condition for positivity under forward Euler integration. The construction accommodates general collision kernels and provides a reusable collision surrogate for both spatially homogeneous and inhomogeneous kinetic solvers, reducing the cost of repeated collision evaluations while retaining these structural properties.

The remainder of this paper is organized as follows. Section~\ref{sec::boltzmann_collision_operator} introduces the Boltzmann collision operator and its numerical formulation. Section~\ref{sec:method} presents the structure-preserving construction and establishes its kinetic properties. Section~\ref{sec:architecture} describes the neural realization and training. Numerical experiments are presented in Section~\ref{sec:experiments}, and concluding remarks are given in Section~\ref{sec:conclusion}.

\section{Boltzmann Collision Operator and Numerical Formulation}
\label{sec::boltzmann_collision_operator}

\subsection{Boltzmann equation and collision operator}
\label{sec::boltzmann}

We consider the nondimensional Boltzmann equation for the nonnegative particle distribution function $f=f(t,\mathbf{x},\mathbf{v})$,
\begin{equation}
  \frac{\partial f}{\partial t} + \sum_{i=1}^d v_i \frac{\partial f}{\partial x_i} = \frac{1}{\mathrm{Kn}}\mathcal{Q}[f,f],
  \label{eq:boltzmann}
\end{equation}
where $t\geq 0$, $\mathbf{v}=(v_1,v_2,v_3)\in\mathbb{R}^{3}$, $\mathbf{x} = (x_1,\ldots,x_d)\in\Omega\subset\mathbb{R}^{d}$ with $d \le 3$, and $\mathrm{Kn}>0$ is the Knudsen number. The collision operator acts pointwise in $t$ and $\mathbf{x}$. Suppressing this dependence, we write
\begin{equation}
  \mathcal{Q}[f,f](\mathbf{v}) = \int_{\mathbb{R}^{3}}\int_{\mathbb{S}^{2}} \mathcal{B}(\mathbf{v}-\mathbf{v}_*,\boldsymbol{\sigma}) \left[f(\mathbf{v}')f(\mathbf{v}_*')-f(\mathbf{v})f(\mathbf{v}_*)\right] \,\mathrm{d}\boldsymbol{\sigma}\,\mathrm{d}\mathbf{v}_*,
  \label{eq:boltzmann_collision}
\end{equation}
where $\mathcal{B}\geq 0$ is the collision kernel. The pairs $(\mathbf{v},\mathbf{v}_*)$ and $(\mathbf{v}',\mathbf{v}_*')$ denote the velocities before and after collision, respectively. For elastic binary collisions, the post-collision velocities are
\begin{equation}
  \mathbf{v}' =\frac{\mathbf{v}+\mathbf{v}_*}{2} +\frac{|\mathbf{v}-\mathbf{v}_*|}{2}\boldsymbol{\sigma},\quad \mathbf{v}_*' =\frac{\mathbf{v}+\mathbf{v}_*}{2} -\frac{|\mathbf{v}-\mathbf{v}_*|}{2}\boldsymbol{\sigma},
  \label{eq:post_collision_velocity}
\end{equation}
where $\boldsymbol{\sigma}\in\mathbb{S}^{2}$ gives the direction of the outgoing relative velocity and $\mathrm{d}\boldsymbol{\sigma}$ denotes the surface measure on $\mathbb{S}^{2}$.

The collision kernel characterizes the intermolecular interaction. A widely used model is the variable hard-sphere (VHS) model, whose kernel is independent of $\boldsymbol{\sigma}$ and has the form
\begin{equation}
  \mathcal{B}(\mathbf{v}-\mathbf{v}_*,\boldsymbol{\sigma}) = C_{\omega}|\mathbf{v}-\mathbf{v}_*|^{2(1-\omega)},
  \label{eq:vhs_kernel}
\end{equation}
where $\omega$ is the viscosity index and $C_{\omega}>0$ is a normalization constant. The choice $\omega=1$ gives a constant collision kernel, known as Maxwell molecules.

The density $\rho$, macroscopic velocity $\mathbf{u}$, and temperature $T$ are defined through the velocity moments of $f$:
\[
\rho=\int_{\mathbb{R}^3} f\,\mathrm{d}\mathbf{v},
\quad
\mathbf{u}=\frac{1}{\rho}\int_{\mathbb{R}^3}\mathbf{v}f\,\mathrm{d}\mathbf{v},
\quad
T=\frac{1}{3\rho}\int_{\mathbb{R}^3}|\mathbf{v}-\mathbf{u}|^2f\,\mathrm{d}\mathbf{v}.
\]

\subsection{Gain and loss decomposition}
\label{sec::gain_loss_structure}

For collision kernels with finite total cross section, the Boltzmann collision operator admits the gain--loss decomposition $\mathcal{Q}[f,f]=\mathcal{Q}^{+}[f]-\mathcal{Q}^{-}[f]$, where the gain term is
\begin{equation}
  \mathcal{Q}^{+}[f](\mathbf{v}) = \int_{\mathbb{R}^{3}}\int_{\mathbb{S}^{2}} \mathcal{B}(\mathbf{v}-\mathbf{v}_*,\boldsymbol{\sigma}) f(\mathbf{v}')f(\mathbf{v}_*') \,\mathrm{d}\boldsymbol{\sigma}\,\mathrm{d}\mathbf{v}_*,
  \label{eq:gain_term}
\end{equation}
and the loss term can be written as $\mathcal{Q}^{-}[f](\mathbf{v})=\nu[f](\mathbf{v})f(\mathbf{v})$, with collision frequency
\begin{equation}
  \nu[f](\mathbf{v}) = \int_{\mathbb{R}^{3}}\int_{\mathbb{S}^{2}} \mathcal{B}(\mathbf{v}-\mathbf{v}_*,\boldsymbol{\sigma}) f(\mathbf{v}_*) \,\mathrm{d}\boldsymbol{\sigma}\,\mathrm{d}\mathbf{v}_*.
  \label{eq:collision_frequency}
\end{equation}
For a nonnegative distribution and a nonnegative collision kernel, both $\mathcal{Q}^{+}[f]$ and $\nu[f]$ are nonnegative. The representation
\begin{equation}
  \mathcal{Q}[f,f](\mathbf{v}) = \mathcal{Q}^{+}[f](\mathbf{v}) - \nu[f](\mathbf{v})f(\mathbf{v})
  \label{eq:gain_loss_frequency}
\end{equation}
separates particle gain from collision loss and exposes a key structural feature of the Boltzmann operator. This decomposition is closely related to the preservation of nonnegativity by the collision dynamics.

For the VHS kernel in~\eqref{eq:vhs_kernel}, the collision frequency reduces to
\[
  \nu[f](\mathbf{v}) = 4\pi C_{\omega} \int_{\mathbb{R}^{3}} |\mathbf{v}-\mathbf{v}_*|^{2(1-\omega)} f(\mathbf{v}_*) \,\mathrm{d}\mathbf{v}_*.
\]
In particular, for the normalized constant kernel $\mathcal{B}=1/(4\pi)$, the collision frequency reduces to $\nu[f](\mathbf{v})=\rho[f]$, where $\rho[f]$ is the density. Thus the collision frequency is independent of $\mathbf{v}$, and the collision operator takes the form $\mathcal{Q}[f,f]=\mathcal{Q}^{+}[f]-\rho[f]f$.

\subsection{Properties of the collision operator}
\label{sec::kinetic_structure}

The Boltzmann collision operator is constrained by several structural properties that determine the macroscopic behavior and long-time evolution of the collision dynamics. We first recall the collision invariants $\phi(\mathbf{v})\in\operatorname{span}\{1,v_1,v_2,v_3,|\mathbf{v}|^2\}$. For every such $\phi$,
\begin{equation}
  \int_{\mathbb{R}^3} \mathcal{Q}[f,f](\mathbf{v})\phi(\mathbf{v})\,\mathrm{d}\mathbf{v}=0.
  \label{eq:collision_invariants}
\end{equation}
The choices $\phi=1$, $\phi=v_i$ for $i=1,2,3$, and $\phi=|\mathbf{v}|^2$ correspond to conservation of mass, momentum, and energy, respectively. Collisions redistribute particles in velocity space while leaving these moments unchanged at each spatial location.

In addition to these conservation laws, the collision dynamics possess a dissipative structure. The Boltzmann entropy is defined by
\begin{equation}
  H(f)=\int_{\mathbb{R}^3} f(\mathbf{v})\log f(\mathbf{v})\,\mathrm{d}\mathbf{v}
  \label{eq:entropy}
\end{equation}
and the $H$-theorem gives
\begin{equation}
  \int_{\mathbb{R}^3} \mathcal{Q}[f,f](\mathbf{v}) \frac{\delta H}{\delta f} \,\mathrm{d}\mathbf{v} =
  \int_{\mathbb{R}^3} \mathcal{Q}[f,f](\mathbf{v})\log f(\mathbf{v})\,\mathrm{d}\mathbf{v}\leq 0.
  \label{eq:h_theorem}
\end{equation}
Hence collisions preserve the macroscopic invariants while driving the distribution in a direction of decreasing $H$. This entropy dissipation gives the collision dynamics an irreversible direction of evolution and describes their relaxation toward equilibrium.

The equilibrium associated with the density $\rho$, macroscopic velocity $\mathbf{u}$, and temperature $T$ is the Maxwellian
\begin{equation}
  \mathcal{M}_{\rho,\mathbf{u},T}(\mathbf{v})=\frac{\rho}{(2\pi T)^{3/2}}\exp\left(-\frac{|\mathbf{v}-\mathbf{u}|^2}{2T}\right),
  \label{eq:maxwellian}
\end{equation}
which satisfies $\mathcal{Q}[\mathcal{M},\mathcal{M}]=0$. Equality in~\eqref{eq:h_theorem} holds only when the distribution function is a Maxwellian. Conservation therefore determines the macroscopic parameters of the equilibrium, while entropy dissipation governs the relaxation toward it. Together with the gain and loss decomposition in the previous subsection, these properties characterize the kinetic structure of the collision dynamics.

\subsection{Fourier spectral approximation}
\label{sec::spectral}

For the numerical evaluation of the collision operator, we use a Fourier spectral method in velocity space. The velocity domain is truncated to $\Omega_L=[-L,L]^3$, and the distribution is extended periodically. To prevent aliasing due to the periodic representation of the collision operator, we truncate the relative velocity at $R$ by replacing $\mathcal{B}(\mathbf{v} - \mathbf{v}_*, \boldsymbol{\sigma})$ with $\mathcal{B}(\mathbf{v} - \mathbf{v}_*, \boldsymbol{\sigma}) \theta(R - |\mathbf{v} - \mathbf{v}_*|)$, where $\theta(\cdot)$ is the Heaviside step function and $R$ satisfies the domain-size condition ${L\geq(3+\sqrt{2})R/4}$~\cite{DimarcoPareschi2014}.

For even $N$, define the grid and mode index sets $\mathcal I_N=\{-N/2,\ldots,N/2-1\}^3$, $\mathcal K_N=\{-N/2,\ldots,N/2\}^3$.
On the grid $\mathbf v_{\boldsymbol\ell}=2L\boldsymbol\ell/N$, $\boldsymbol\ell\in\mathcal I_N$, the trigonometric interpolant is
\begin{equation}
  f_N(\mathbf v)=\frac{1}{(2L)^3}
  \sum_{\mathbf k\in\mathcal K_N}
  \widehat f_{\mathbf k}
  e^{\mathrm i\pi\mathbf k\cdot\mathbf v/L},
  \label{eq:fourier_expansion}
\end{equation}
with coefficients
\begin{equation}
  \widehat f_{\mathbf k}
  =\frac{(2L/N)^3}{c_{\mathbf k}}
  \sum_{\boldsymbol\ell\in\mathcal I_N}
  f(\mathbf v_{\boldsymbol\ell})
  e^{-2\pi\mathrm i\mathbf k\cdot\boldsymbol\ell/N}.
  \label{eq:discrete-fourier-coefficients}
\end{equation}
Here $c_{\mathbf k}=2^{\#\{j:|k_j|=N/2\}}$ accounts for the repeated Nyquist modes.

The fast spectral method exploits the convolution structure of the collision operator after introducing the relative velocity and expressing it in spherical coordinates. The radial integral is discretized using a Jacobi--Gauss--Radau quadrature~\cite{ShenTangWang2011}, while the angular integration is approximated using Lebedev quadrature~\cite{Lebedev1976}:
\begin{equation}
  \int_0^R g^2\psi(g)\,\mathrm{d}g\approx\sum_{j=1}^{J}w_j^g\psi(g_j),\quad \int_{\mathbb{S}^2}\varphi(\boldsymbol{\sigma})\,\mathrm{d}\boldsymbol{\sigma}\approx\sum_{m=1}^{M}w_m^\sigma\varphi(\boldsymbol{\sigma}_m).
  \label{eq:spectral_quadrature}
\end{equation}
Here $g$ is the relative speed, $J$ and $M$ are the numbers of radial and angular nodes, and the weights include the corresponding integration measures. For fixed radial and angular quadrature nodes, the resulting sums over Fourier modes have convolution form and can be evaluated efficiently using fast Fourier transforms~\cite{fast_spectral_method_boltzmann1,fast_spectral_method_boltzmann2}. The overall computational complexity is $\mathcal{O}(JMN^3\log N)$, which reduces to $\mathcal{O}(MN^4\log N)$ when $J=\mathcal{O}(N)$. Despite its high accuracy for smooth distributions, the standard Fourier spectral method does not preserve all the kinetic structure of the collision operator. It conserves mass, but momentum and energy are generally not conserved exactly without additional corrections~\cite{PareschiRey2022}. Positivity and a discrete $H$-theorem are likewise not guaranteed~\cite{CaiFanYing2018}.

The fast spectral method provides reference values of \(\mathcal{Q}[f,f]\) for numerical comparisons. The corresponding velocity discretization and spectral parameters are specified together with the numerical experiments.

\subsection[BGK relaxation and moment-matched equilibrium]{BGK relaxation and moment-matched equilibrium}
\label{sec:bgk-background}
\leavevmode\par\nobreak

The BGK model~\cite{BGK1} approximates the Boltzmann collision operator by relaxation toward a local Maxwellian,
\begin{equation}
  \mathcal{Q}_{\mathrm{BGK}}[f]=\nu_{\mathrm{BGK}}[f](\mathcal{M}[f]-f),
\end{equation}
where the positive relaxation frequency is independent of velocity. On the truncated velocity domain $\Omega_L$, the full-space Maxwellian does not generally match the moments of $f$ evaluated over $\Omega_L$. To ensure moment conservation, we instead determine the equilibrium by entropy minimization.

Let $\rho[f]=\int_{\Omega_L}f\,\mathrm{d}\mathbf{v}$ and define the scaled collision-invariant vector
\begin{equation}
  \psi(\mathbf{v})=\left(1,\frac{v_1}{L},\frac{v_2}{L},\frac{v_3}{L},\frac{|\mathbf{v}|^2}{L^2}\right)^{\!T}.
  \label{eq:collision-invariants}
\end{equation}
For any nonnegative $f\in L^1(\Omega_L)$ with $\rho[f]>0$, define the moment-matched equilibrium by
\begin{equation}
  M[f]=\underset{\substack{h\geq0\\ \int_{\Omega_L}\psi h\,\mathrm{d}\mathbf{v}=\int_{\Omega_L}\psi f\,\mathrm{d}\mathbf{v}}}{\operatorname{argmin}}\int_{\Omega_L}h(\log h-1)\,\mathrm{d}\mathbf{v}.
  \label{eq:moment-matched-equilibrium}
\end{equation}
The admissible densities are integrable with finite entropy, and we use $0\log0=0$. The minimizer exists uniquely and has the exponential form $M[f]=e^{\alpha[f]^T\psi}$, with $\alpha[f]\in\mathbb{R}^5$ uniquely determined by the moment constraints~\cite{discrete_bgk_model2}.

\section{Structure-preserving collision operator}
\label{sec:method}
We first formulate a general structure-preserving construction for approximating the Boltzmann collision operator, independent of how the underlying approximation is parameterized. The construction combines a gain--loss proposal with a conservative moment projection and an entropy correction toward a BGK reference. The proposal preserves equilibrium balance and density scaling, while the two corrections enforce conservation and entropy dissipation. We establish the structural properties and relaxation behavior of this construction, and Section~\ref{sec:architecture} gives a neural realization of the admissible proposal class.

All moments, entropy integrals, and projections below are evaluated on $\Omega_L$, using the invariant vector and equilibrium defined in Section~\ref{sec:bgk-background}.
The state space is
\begin{equation}
  \mathcal{A}_+=\left\{f\in L^\infty(\Omega_L):f\geq0,\ \rho[f]>0\right\}.
  \label{eq:nonnegative-class}
\end{equation}

\subsection[Construction of the structure-preserving operator]{Construction of the structure-preserving operator}
\label{sec:operator-construction}
\label{sec:gain-loss-proposal}
\label{sec:moment-projection}
\label{sec:entropy-correction}

We begin with maps $G,\nu:\mathcal{A}_+\to L^\infty(\Omega_L)$ and the proposal
\begin{equation}
  q[f]=G[f]-\nu[f]f.
  \label{eq:raw-gain-loss}
\end{equation}
For each $f\in\mathcal{A}_+$, we require $\operatorname*{ess\,inf}_{\Omega_L}G[f]>0$ and $\nu[f]>0$ almost everywhere. We also require equilibrium balance,
\begin{equation}
  G[M[f]]=\nu[M[f]]M[f],
  \label{eq:proposal-balance}
\end{equation}
and the physical density scaling,
\begin{equation}
  G[af]=a^2G[f],\quad \nu[af]=a\nu[f],\quad a>0.
  \label{eq:proposal-scaling}
\end{equation}
These conditions give $q[M[f]]=0$ and $q[af]=a^2q[f]$. We next impose conservation and entropy dissipation while retaining a gain--loss representation.

To enforce moment conservation, for $f\in\mathcal{A}_+$ define the weighted projection
\begin{equation}
  H_f=\int_{\Omega_L}f\psi\psi^T\,\mathrm{d}\mathbf{v},\quad \mathcal{P}_fq=q-\left(\psi^TH_f^{-1}\int_{\Omega_L}\psi q\,\mathrm{d}\mathbf{v}\right)f.
  \label{eq:continuous-projection}
\end{equation}
The corrected proposal is $q^{\mathrm c}[f]=\mathcal{P}_fq[f]$. Since the correction is a bounded coefficient multiplied by $f$, it modifies only the loss coefficient and leaves the gain unchanged. This form will be used in the positivity estimates.

\begin{proposition}[Weighted conservative projection]
\label{prop:conservative-projection}
For $f\in\mathcal{A}_+$, $H_f$ is positive definite and $\int_{\Omega_L}\psi\mathcal{P}_fq\,\mathrm{d}\mathbf{v}=0$ for every $q\in L^2(\Omega_L)$. Moreover, $\mathcal{P}_{af}(a^2q)=a^2\mathcal{P}_fq$ for every $a>0$.
\end{proposition}
\begin{proof}
For nonzero $b\in\mathbb{R}^5$, the polynomial $b^T\psi$ cannot vanish on a set of positive volume. Since $f\geq0$ has positive mass,
\[
  b^TH_fb=\int_{\Omega_L}f(b^T\psi)^2\,\mathrm{d}\mathbf{v}>0.
\]
Taking moments in~\eqref{eq:continuous-projection} gives
\[
  \int_{\Omega_L}\psi\mathcal{P}_fq\,\mathrm{d}\mathbf{v}
  =\int_{\Omega_L}\psi q\,\mathrm{d}\mathbf{v}-H_fH_f^{-1}\int_{\Omega_L}\psi q\,\mathrm{d}\mathbf{v}=0.
\]
Finally, $H_{af}=aH_f$, so direct substitution gives the scaling identity.
\end{proof}

The moment correction does not enforce entropy dissipation. To impose an entropy bound without changing the moments, we use the moment-matched BGK reference
\begin{equation}
  Q_{\mathrm B}[f]=r_{\mathrm B}\rho[f]T[f]^{1-\omega}(M[f]-f),\quad r_{\mathrm B}>0.
  \label{eq:bgk-reference}
\end{equation}
Here $r_{\mathrm B}$ is constant for the prescribed VHS kernel, and $T[f]>0$ is the temperature computed from the moments on $\Omega_L$. Uniqueness of the moment-matched exponential gives $M[af]=aM[f]$ for $a>0$. Since $T[af]=T[f]$, the density factor and the scaling of $M[f]-f$ give $Q_{\mathrm B}[af]=a^2Q_{\mathrm B}[f]$.

To compare the entropy dissipation of the corrected proposal and the BGK reference, we use the entropy pairing $D_f(q)$. The relative entropy $\mathscr{H}(f\mid M[f])$ measures departure from equilibrium. For $f\in\mathcal{A}_+$ and a collision direction $q$ whose five collision-invariant moments vanish, define
\begin{equation}
  D_f(q)=\int_{\Omega_L}q\log\frac{f}{\rho[f]}\,\mathrm{d}\mathbf{v},\quad \mathscr{H}(f\mid M[f])=\int_{\Omega_L}f\log\frac{f}{M[f]}\,\mathrm{d}\mathbf{v}.
  \label{eq:entropy-pairing}
\end{equation}
We use $0\log0=0$ in the relative entropy. At zeros of $f$, the pairing integrand is defined as $-\infty$ if $q>0$ and as zero if $q=0$.

Since the five collision-invariant moments of $q$ vanish and $\log M[f]=\alpha[f]^T\psi$,
\begin{equation}
  D_f(q)=\int_{\Omega_L}q\log f\,\mathrm{d}\mathbf{v}=\int_{\Omega_L}q\log\frac{f}{M[f]}\,\mathrm{d}\mathbf{v}.
  \label{eq:relative-entropy-pairing}
\end{equation}
The identities include the case $D_f(q)=-\infty$. Also, $0\leq\mathscr{H}(f\mid M[f])<\infty$, with equality only at $f=M[f]$.

For the two collision terms considered here, these pairings are well defined when $f$ vanishes. The projection leaves the gain of $q^{\mathrm c}[f]$ equal to $G[f]$, while the BGK gain is $r_{\mathrm B}\rho[f]T[f]^{1-\omega}M[f]$. Both are bounded above and away from zero, by the proposal assumption and the exponential form of $M[f]$, respectively. Their loss coefficients are bounded, so the loss contributions to the pairing are finite because $f\log(f/\rho[f])$ is integrable. The gain contributions are bounded above and are finite exactly when $\log f\in L^1(\Omega_L)$. Thus both pairings belong to $[-\infty,\infty)$ and are finite under the same condition, which ensures that the mixing weight below is well defined.

\begin{lemma}[Dissipation of the BGK reference]
\label{lem:bgk-dissipation}
For $f\in\mathcal{A}_+$,
\begin{equation}
  \begin{aligned}
    D_f(Q_{\mathrm B}[f])&=-r_{\mathrm B}\rho[f]T[f]^{1-\omega}\bigl[\mathscr{H}(f\mid M[f])+\mathscr{H}(M[f]\mid f)\bigr]\\
    &\leq-r_{\mathrm B}\rho[f]T[f]^{1-\omega}\mathscr{H}(f\mid M[f]).
  \end{aligned}
  \label{eq:bgk-entropy-identity}
\end{equation}
The identity allows $\mathscr{H}(M[f]\mid f)=+\infty$. In particular, $D_f(Q_{\mathrm B}[f])<0$ whenever $f\ne M[f]$.
\end{lemma}
\begin{proof}
Using~\eqref{eq:relative-entropy-pairing},
\[
  D_f(Q_{\mathrm B}[f])=-r_{\mathrm B}\rho[f]T[f]^{1-\omega}\int_{\Omega_L}(f-M[f])\log\frac{f}{M[f]}\,\mathrm{d}\mathbf{v}.
\]
Expanding into the two relative entropies and using their nonnegativity proves the result.
\end{proof}

Fix $\kappa\in[0,1]$ and set the entropy target
\begin{equation}
  T_\kappa[f]=-\kappa r_{\mathrm B}\rho[f]T[f]^{1-\omega}\mathscr{H}(f\mid M[f]).
  \label{eq:entropy-target}
\end{equation}
Lemma~\ref{lem:bgk-dissipation} gives $D_f(Q_{\mathrm B}[f])\leq T_\kappa[f]$, so the BGK reference always meets the target. We mix it with the corrected proposal only when the latter exceeds this target. Set
\begin{equation}
  \gamma_\kappa[f]=\begin{cases}0,&D_f(q^{\mathrm c}[f])\leq T_\kappa[f],\\[4pt]\displaystyle\frac{D_f(q^{\mathrm c}[f])-T_\kappa[f]}{D_f(q^{\mathrm c}[f])-D_f(Q_{\mathrm B}[f])},&D_f(q^{\mathrm c}[f])>T_\kappa[f].\end{cases}
  \label{eq:entropy-gate}
\end{equation}
The corrected operator is
\begin{equation}
  \mathcal{Q}_\kappa[f]=(1-\gamma_\kappa[f])q^{\mathrm c}[f]+\gamma_\kappa[f]Q_{\mathrm B}[f].
  \label{eq:final-gated-operator}
\end{equation}
A weight in $[0,1]$ preserves the vanishing collision-invariant moments and nonnegative gain of both terms. The choice $\kappa=0$ enforces non-strict entropy dissipation, while $\kappa>0$ imposes a relative entropy gap that will exclude nonequilibrium stationary states.

\begin{proposition}[Entropy correction]
\label{prop:entropy-correction}
For every $f\in\mathcal{A}_+$ and $\kappa\in[0,1]$, the weight in~\eqref{eq:entropy-gate} is well defined and belongs to $[0,1]$. The corrected operator satisfies $D_f(\mathcal{Q}_\kappa[f])\leq T_\kappa[f]$.
\end{proposition}
\begin{proof}
If $D_f(q^{\mathrm c}[f])\leq T_\kappa[f]$, the choice $\gamma_\kappa[f]=0$ gives the result. Otherwise both pairings are finite by the criterion above, and Lemma~\ref{lem:bgk-dissipation} gives
\[
  D_f(Q_{\mathrm B}[f])\leq T_\kappa[f]<D_f(q^{\mathrm c}[f]).
\]
Hence $0<\gamma_\kappa[f]\leq1$, and linearity gives
\[
  D_f(\mathcal{Q}_\kappa[f]) =(1-\gamma_\kappa[f])D_f(q^{\mathrm c}[f])+\gamma_\kappa[f]D_f(Q_{\mathrm B}[f]) =T_\kappa[f].
\]
\end{proof}

\begin{remark}[Activity of the entropy correction]
In our finite dataset, every sample satisfies $D_f(q^{\mathrm c}[f])<0$. With $\kappa=0$, the weight is therefore $\gamma_0[f]=0$, and the entropy correction leaves the corrected proposal unchanged. We retain the correction to guarantee the entropy inequality for every admissible input, independently of the learned parameters.
\end{remark}

Both corrections retain the gain--loss structure. Indeed,
\begin{equation}
\begin{aligned}
  \mathcal{Q}_\kappa[f]
  ={}&(1-\gamma_\kappa[f])G[f]+\gamma_\kappa[f]r_{\mathrm B}\rho[f]T[f]^{1-\omega}M[f]\\
  &-\biggl[(1-\gamma_\kappa[f])
  \left(\nu[f]+\psi^TH_f^{-1}\int_{\Omega_L}\psi q[f]\,\mathrm{d}\mathbf{v}\right)
  +\gamma_\kappa[f]r_{\mathrm B}\rho[f]T[f]^{1-\omega}\biggr]f.
\end{aligned}
\label{eq:corrected-gain-loss}
\end{equation}
The gain remains positive, while the bounded coefficient multiplying $f$ may have either sign because of the moment projection.

\subsection[Structural properties and positivity]{Structural properties and positivity}
\label{sec:operator-properties}
\label{sec:dynamics}

We now establish the kinetic structure of the corrected operator, derive positivity bounds, and prove exponential relaxation for $\kappa>0$.

\begin{theorem}[Exact kinetic structure]
\label{prop:structural-identities}
For $f\in\mathcal{A}_+$, $\kappa\in[0,1]$, and any admissible proposal from Section~\ref{sec:gain-loss-proposal}, the operator~\eqref{eq:final-gated-operator} satisfies
\begin{align}
  \int_{\Omega_L}\psi\mathcal{Q}_\kappa[f]\,\mathrm{d}\mathbf{v}&=0,\quad \mathcal{Q}_\kappa[M[f]]=0, \label{eq:operator-conservation}\\
  D_f(\mathcal{Q}_\kappa[f])&\leq-\kappa r_{\mathrm B}\rho[f]T[f]^{1-\omega}\mathscr{H}(f\mid M[f]), \label{eq:operator-entropy}\\
  \mathcal{Q}_\kappa[af]&=a^2\mathcal{Q}_\kappa[f],\quad a>0. \label{eq:operator-scaling}
\end{align}
If $\kappa>0$, then $\mathcal{Q}_\kappa[f]=0$ if and only if $f=M[f]$.
\end{theorem}
\begin{proof}
Since $M[f]$ has the same moments as $f$, the two inputs define the same entropy-minimization problem. Uniqueness of its minimizer gives $M[M[f]]=M[f]$. Both $q^{\mathrm c}[f]$ and $Q_{\mathrm B}[f]$ have vanishing collision-invariant moments and vanish at equilibrium. Their convex combination therefore satisfies~\eqref{eq:operator-conservation}. The entropy bound follows from Proposition~\ref{prop:entropy-correction}.

Under $f\mapsto af$, both $\rho[f]$ and $M[f]$ scale by $a$, while the projected proposal and the BGK reference scale by $a^2$. Moreover,
\[
  D_{af}(a^2q)=a^2D_f(q),\quad \mathscr{H}(af\mid M[af])=a\mathscr{H}(f\mid M[f]),\quad T_\kappa[af]=a^2T_\kappa[f].
\]
Thus $\gamma_\kappa[af]=\gamma_\kappa[f]$, which proves~\eqref{eq:operator-scaling}. If $\kappa>0$ and $\mathcal{Q}_\kappa[f]=0$, the entropy bound forces $\mathscr{H}(f\mid M[f])=0$, hence $f=M[f]$. The converse follows from equilibrium preservation.
\end{proof}

Denote the gain and the bracketed coefficient in~\eqref{eq:corrected-gain-loss} by $G_{\mathrm{eff}}[f]$ and $\Lambda_\kappa[f]$, respectively. Then
\begin{equation}
  \mathcal{Q}_\kappa[f]=G_{\mathrm{eff}}[f]-\Lambda_\kappa[f]f.
  \label{eq:effective-gain-loss}
\end{equation}
The following positivity bounds use only $G_{\mathrm{eff}}[f]\geq0$ and an upper bound on $\Lambda_\kappa[f]$.

\begin{proposition}[Forward Euler positivity]
\label{prop:euler-positivity}
Let $f\in\mathcal{A}_+$ and define
\[
  b[f]=\|(\Lambda_\kappa[f])_+\|_{L^\infty},\quad x_+=\max\{x,0\}.
\]
If $\Delta t\geq0$ and $\Delta t\,b[f]\leq1$, then the collision update $f^+=f+\Delta t\,\mathcal{Q}_\kappa[f]$ belongs to $\mathcal{A}_+$ and has the same five moments as $f$.
\end{proposition}
\begin{proof}
The effective gain--loss representation gives
\[
  f^+=(1-\Delta t\,\Lambda_\kappa[f])f+\Delta t\,G_{\mathrm{eff}}[f]\geq(1-\Delta t\,b[f])f+\Delta t\,G_{\mathrm{eff}}[f]\geq0.
\]
The update is bounded and preserves the five moments by Theorem~\ref{prop:structural-identities}, so $f^+\in\mathcal{A}_+$.
\end{proof}

For the scaled collision equation $\partial_t f=\mathcal{Q}_\kappa[f]/\mathrm{Kn}$, the sufficient condition is $\Delta t\,b[f]\leq\mathrm{Kn}$.

In the trajectory results below, a classical solution on $[0,t_*)$ satisfies
\[
  f\in C^1([0,T];L^\infty(\Omega_L))\quad\text{for every }0<T<t_*.
\]
We assume its effective gain and loss coefficient are strongly measurable with locally integrable $L^\infty$ norms.

\begin{proposition}[Positivity along collision trajectories]
\label{prop:continuous-dynamics}
Let $f(t)$ be a classical solution of $\partial_t f=\mathcal{Q}_\kappa[f]$ in $\mathcal{A}_+$ on $[0,t_*)$, with initial state $f_0$. Set $b(t)=\|(\Lambda_\kappa[f(t)])_+\|_\infty$. For $0<t<t_*$,
\begin{equation}
  f(t,\mathbf{v})\geq e^{-\int_0^t b(s)\,\mathrm{d}s}f_0(\mathbf{v})\geq0
  \label{eq:positive-trajectory-bound}
\end{equation}
almost everywhere. If $t_*<\infty$, $\int_0^{t_*}b(s)\,\mathrm{d}s<\infty$, and an $L^\infty$ limit exists, the same bound holds at $t=t_*$.
\end{proposition}
\begin{proof}
The effective gain--loss representation and $b(t)\geq\Lambda_\kappa[f(t)]$ give
\[
  \partial_t f+b(t)f=G_{\mathrm{eff}}[f(t)]+\bigl(b(t)-\Lambda_\kappa[f(t)]\bigr)f\geq0.
\]
Hence
\[
  \frac{\mathrm{d}}{\mathrm{d}t}\left[e^{\int_0^t b(s)\,\mathrm{d}s}f(t,\mathbf{v})\right]\geq0.
\]
Integration gives~\eqref{eq:positive-trajectory-bound}, and passage to the $L^\infty$ limit gives the endpoint statement.
\end{proof}

An exponential Euler update preserves nonnegativity for any $\Delta t>0$. Rewrite the collision equation as
\[
  \partial_t f+b(t)f=\mathcal{Q}_\kappa[f]+b(t)f.
\]
Freezing $b(t)$ at $b^n=b[f^n]$ and the right-hand side at $\mathcal{Q}_\kappa[f^n]+b^nf^n$, then integrating over one time step, gives
\begin{equation}
  f^{n+1}=f^n+\alpha_n\mathcal{Q}_\kappa[f^n],\quad \alpha_n=\frac{1-e^{-b^n\Delta t}}{b^n},
  \label{eq:exponential-positivity-update}
\end{equation}
where $\alpha_n=\Delta t$ when $b^n=0$. Since $\alpha_n\geq0$ and $\alpha_n b^n\leq1$, Proposition~\ref{prop:euler-positivity} gives $f^{n+1}\in\mathcal{A}_+$ and exact preservation of the five moments. Moreover, $\alpha_n=\Delta t+\mathcal{O}(\Delta t^2)$, so the update is first-order consistent.

\medskip
\noindent\textbf{Homogeneous relaxation.}

For $\kappa>0$, conservation and the entropy gap give exponential relaxation to the equilibrium determined by the initial moments.

\begin{theorem}[Exponential relaxation]
\label{thm:homogeneous-relaxation}
Fix $0<\kappa\leq1$. Let $f(t)$ be a global classical solution of $\partial_t f=\mathcal{Q}_\kappa[f]$ in $\mathcal{A}_+$ with initial state $f_0$. Set $\rho_0=\rho[f_0]$, $T_0=T[f_0]$, and $M_0=M[f_0]$. Then $\rho[f(t)]=\rho_0$, $T[f(t)]=T_0$, $M[f(t)]=M_0$, and
\begin{align}
  \mathscr{H}(f(t)\mid M_0)&\leq e^{-\kappa r_{\mathrm B}\rho_0T_0^{1-\omega}t}\mathscr{H}(f_0\mid M_0), \label{eq:relative-entropy-decay}\\
  \|f(t)-M_0\|_{L^1(\Omega_L)}&\leq\sqrt{2\rho_0\mathscr{H}(f_0\mid M_0)}\,e^{-\kappa r_{\mathrm B}\rho_0T_0^{1-\omega}t/2}. \label{eq:l1-relaxation}
\end{align}
In particular, $f(t)\to M_0$ in $L^1(\Omega_L)$ as $t\to\infty$.
\end{theorem}
\begin{proof}
The five conserved moments give $\rho[f(t)]=\rho_0$, $T[f(t)]=T_0$, and $M[f(t)]=M_0$. Set $H(t)=\mathscr{H}(f(t)\mid M_0)$. Boundedness of $f$ and the positive lower bound on $M_0$ ensure that $H$ is finite. It is also continuous, including at zeros of $f$, since $f$ is continuous in $L^\infty$ and $0\log0=0$. We show that
\[
  H(t)-H(s)\leq-\kappa r_{\mathrm B}\rho_0T_0^{1-\omega}\int_s^t H(\tau)\,\mathrm{d}\tau,\quad 0\leq s<t.
\]
To justify this inequality at zeros of $f$, set $H_\varepsilon(t)=\int_{\Omega_L}(f(t)+\varepsilon)\log((f(t)+\varepsilon)/M_0)\,\mathrm{d}\mathbf{v}$. Mass conservation yields
\[
  H_\varepsilon'(t)=\int_{\Omega_L}\mathcal{Q}_\kappa[f(t)]\log\frac{f(t)+\varepsilon}{M_0}\,\mathrm{d}\mathbf{v}.
\]
Denote this integrand by $I_\varepsilon(t,\mathbf{v})$. On $[s,t]$, the gain--loss form gives the uniform upper bound
\[
  I_\varepsilon(\tau,\mathbf{v})\leq B(\tau):=C\bigl(\|G_{\mathrm{eff}}[f(\tau)]\|_\infty+\|\Lambda_\kappa[f(\tau)]\|_\infty\bigr),\quad 0<\varepsilon\leq1,
\]
with $B\in L^1(s,t)$. Apply Fatou's lemma to the nonnegative functions $B-I_\varepsilon$ on $[s,t]\times\Omega_L$. Subtracting the finite integral of $B$ gives
\[
  \limsup_{\varepsilon\downarrow0}\int_s^t\!\int_{\Omega_L}I_\varepsilon\,\mathrm{d}\mathbf{v}\,\mathrm{d}\tau
  \leq\int_s^t\!\int_{\Omega_L}\lim_{\varepsilon\downarrow0}I_\varepsilon\,\mathrm{d}\mathbf{v}\,\mathrm{d}\tau.
\]
The left-hand side equals $H(t)-H(s)$, since $H_\varepsilon\to H$ at both endpoints. The inner integral on the right is the extended entropy pairing $D_{f(\tau)}(\mathcal{Q}_\kappa[f(\tau)])$ by~\eqref{eq:relative-entropy-pairing}. Hence
\[
  H(t)-H(s)\leq\int_s^t D_{f(\tau)}\bigl(\mathcal{Q}_\kappa[f(\tau)]\bigr)\,\mathrm{d}\tau
  \leq-\kappa r_{\mathrm B}\rho_0T_0^{1-\omega}\int_s^t H(\tau)\,\mathrm{d}\tau,
\]
where the last step uses Theorem~\ref{prop:structural-identities}. Gronwall's inequality proves~\eqref{eq:relative-entropy-decay}. Finally, the equal-mass Pinsker inequality
\[
  \|f(t)-M_0\|_{L^1}^2\leq2\rho_0 H(t)
\]
gives~\eqref{eq:l1-relaxation} and the asserted convergence.
\end{proof}

For $\kappa=0$, the same argument gives nonincreasing relative entropy. For the scaled equation $\partial_t f=\mathcal{Q}_\kappa[f]/\mathrm{Kn}$, the rates in~\eqref{eq:relative-entropy-decay}--\eqref{eq:l1-relaxation} are divided by $\mathrm{Kn}$.

\section{Neural realization and training}
\label{sec:architecture}

The construction in Section~\ref{sec:method} gives sufficient conditions for a learned collision operator to preserve the desired kinetic structure. We parameterize the gain term and the collision frequency by neural networks with parameters $\theta$, writing $G_\theta$, $\nu_\theta$, and $q_\theta$ for the resulting maps and $\mathcal{Q}_{\theta,\kappa}$ for the corrected operator. The parameterization satisfies the admissibility conditions by design. The neural network learns the dependence on the normalized velocity distribution, while the density scaling and structural constraints are enforced explicitly. We then describe the discrete realization and training objective.

\subsection{State normalization and encoding}
\label{sec:state-encoding}

To obtain the homogeneity in~\eqref{eq:proposal-scaling}, we construct a density-invariant encoding that preserves the distribution shape, its departure from equilibrium, and the local velocity-space structure. Let $\mathbf{u}$ and $T$ be the macroscopic velocity and temperature of $f$ on $\Omega_L$. Set $\mathbf{c}(\mathbf{v})=(\mathbf{v}-\mathbf{u})/\sqrt{T}$ and $A[f]=\|\max\{f,M[f]\}\|_{L^\infty(\Omega_L)}$. With the small constant $\varepsilon_{\mathrm f}=10^{-12}$, define
\begin{equation}
  \mathcal{E}[f](\mathbf{v})=\left(\frac{f}{A[f]},\frac{M[f]}{A[f]},\frac{f-M[f]}{f+M[f]+\varepsilon_{\mathrm f}A[f]},\frac{\mathbf{c}}{\sqrt{1+|\mathbf{c}|^2}},\frac{|\mathbf{c}|^2}{1+|\mathbf{c}|^2}\right)^{\!T}.
  \label{eq:seven-channel-encoding}
\end{equation}
The seven channels capture the normalized distribution, its relative deviation from equilibrium, and velocity coordinates in local thermal units. The term $\varepsilon_{\mathrm f}A[f]$ stabilizes the relative difference where both $f$ and $M[f]$ are small compared with $A[f]$. All components are bounded. Since $M[af]=aM[f]$ and $A[af]=aA[f]$, while $\mathbf{u}$ and $T$ are unchanged, $\mathcal{E}[af]=\mathcal{E}[f]$ for every $a>0$.

\subsection{Fourier neural operator parameterization}
\label{sec:fno-backbone}

The gain and collision frequency at a given velocity depend nonlocally on the distribution over velocity space. We represent this dependence using a Fourier neural operator (FNO)~\cite{FNO}, a neural-operator architecture that couples velocity locations through spectral convolution~\cite{Kovachki2023NeuralOperator}. Since the collision operator is local in physical space, the same velocity-space mapping is applied independently at each spatial location.

The encoded state is first lifted pointwise to $w$ latent channels,
\[
  z_0(\mathbf{v})=P\mathcal{E}[f](\mathbf{v})+p,
\]

where $P\in\mathbb{R}^{w\times7}$ and $p\in\mathbb{R}^w$. The latent representation is then propagated through $K\geq1$ Fourier layers,
\begin{equation}
  z_{\ell+1}(\mathbf{v})=\phi_\ell\!\left(W_\ell z_\ell(\mathbf{v})+b_\ell+\mathcal{F}^{-1}\!\left[R_\ell\mathcal{F}z_\ell\right](\mathbf{v})\right),\quad 0\leq\ell<K.
  \label{eq:fourier-layer}
\end{equation}
Here $\mathcal{F}$ denotes the Fourier-series transform with respect to the velocity variable on $\Omega_L$, applied separately to each of the $w$ feature channels. The matrices $W_\ell\in\mathbb{R}^{w\times w}$ and vectors $b_\ell\in\mathbb{R}^w$ define pointwise affine maps. At each retained Fourier mode $\mathbf{k}$, the learned matrix $R_\ell(\mathbf{k})\in\mathbb{C}^{w\times w}$ mixes the transformed channels. This multiplier is supported on finitely many low modes. We use $\phi_\ell=\operatorname{GELU}$~\cite{GELU} for $0\leq\ell<K-1$ and the identity map for $\phi_{K-1}$.

Two pointwise affine readouts of the shared latent field give
\begin{equation}
  a_\theta[f](\mathbf{v})=w_g^Tz_K(\mathbf{v})+\beta_g,\quad b_\theta[f](\mathbf{v})=w_r^Tz_K(\mathbf{v})+\beta_r.
  \label{eq:dual-head-operator}
\end{equation}
All parameters in the lifting, Fourier layers, and readouts are trainable and collectively denoted by $\theta$.
Since $\mathcal{E}[af]=\mathcal{E}[f]$, both output fields are invariant under positive density scaling. These unconstrained outputs are converted into the admissible gain and collision frequency in the next subsection.

\subsection{Gain and collision-frequency parameterization}
\label{sec:gain-frequency-parameterization}

We enforce equilibrium balance by making the learned perturbations vanish when $f=M[f]$. Define
\begin{equation}
  g[f]=\frac{\|f-M[f]\|_{L^2(\Omega_L)}}{\|f\|_{L^2(\Omega_L)}+\|M[f]\|_{L^2(\Omega_L)}}.
  \label{eq:equilibrium-gate}
\end{equation}
Then $g[M[f]]=0$ and $g[af]=g[f]$ for $a>0$. This function will play the role of suppressing the learned deviations at equilibrium without introducing additional density dependence.

Let $\sigma(s)=(1+e^{-s})^{-1}$ and $\operatorname{logit}(p)=\log(p/(1-p))$. With the constant $r_{\mathrm B}$ from~\eqref{eq:bgk-reference}, choose $A_g,A_r>0$ and $0<r_{\min}<r_{\mathrm B}<r_{\max}$, and set $\ell_0=\operatorname{logit}((r_{\mathrm B}-r_{\min})/(r_{\max}-r_{\min}))$.

We define the learned collision frequency by
\begin{equation}
  \nu_\theta[f]=\rho[f]T[f]^{1-\omega}\left[r_{\min}+(r_{\max}-r_{\min})\sigma\left(\ell_0+A_r\tanh\frac{g[f]b_\theta[f]}{A_r}\right)\right].
  \label{eq:learned-frequency}
\end{equation}
The logistic readout bounds $\nu_\theta[f]/(\rho[f]T[f]^{1-\omega})$ between $r_{\min}$ and $r_{\max}$, while $A_r$ controls the logit perturbation. When $b_\theta[f]=0$, substituting the definition of $\ell_0$ gives the BGK reference frequency $r_{\mathrm B}\rho[f]T[f]^{1-\omega}$.

For the gain term $G_{\theta}[f]$, we adopt the form
\begin{equation}
  G_\theta[f]=C[f]\,\sigma\left(\operatorname{logit}\frac{\nu_\theta[f]M[f]}{C[f]}+A_g\tanh\frac{g[f]a_\theta[f]}{A_g}\right).
  \label{eq:bounded-gain}
\end{equation}
When $f = M[f]$, the function $g[f]$ vanishes, yielding $G_{\theta}[f] = \nu_{\theta}[f] M[f]$, identical the gain term in the BGK operator. Therefore, the parameter $A_g$ controls the deviation from the BGK operator owing to the non-equilibrium part of $f$. The function $C[f]$, depending on $\mathbf{c}$, is an envelope function that avoids excessive values in the gain term. For $\tau>0$, we choose
\begin{equation}
  C[f]=r_{\max}\rho[f]T[f]^{1-\omega}\left[\max\{f,M[f]\}+e^{-\tau|\mathbf{c}|^2}\bigl\|\max\{f,M[f]\}e^{\tau|\mathbf{c}|^2}\bigr\|_{L^\infty(\Omega_L)}\right].
  \label{eq:gain-envelope}
\end{equation}
Since $r_{\max}>r_{\mathrm B}$, this envelope also bounds the BGK gain $r_{\mathrm B}\rho[f]T[f]^{1-\omega}M[f]$.  Such an envelope also controls the velocity tail of the gain term, avoiding growing tails when $f$ decays.

The following proposition verifies that $q_\theta[f]=G_\theta[f]-\nu_\theta[f]f$ belongs to the admissible class of Section~\ref{sec:method}.

\begin{proposition}[Admissibility of the neural parameterization]
\label{prop:neural-realization}
For every $f\in\mathcal{A}_+$ and every parameter vector $\theta$ with finite entries, the preceding parameterization produces bounded fields $G_\theta[f]$ and $\nu_\theta[f]$ that are bounded away from zero and satisfy~\eqref{eq:proposal-balance}--\eqref{eq:proposal-scaling}. If both network outputs vanish, $a_\theta[f]=b_\theta[f]=0$, then $q_\theta[f]=Q_{\mathrm B}[f]$.
\end{proposition}
\begin{proof}
For $f\in\mathcal{A}_+$, positive mass cannot be concentrated at a single velocity, so $T>0$. The exponential form of $M[f]$ ensures that it is bounded above and away from zero. The bounded encoding and finite Fourier layers give bounded output fields $a_\theta[f]$ and $b_\theta[f]$. The frequency satisfies $r_{\min}\rho[f]T[f]^{1-\omega}<\nu_\theta[f]<r_{\max}\rho[f]T[f]^{1-\omega}$, and the envelope satisfies
\[
  C[f]\geq2r_{\max}\rho[f]T[f]^{1-\omega}\max\{f,M[f]\}\geq2r_{\max}\rho[f]T[f]^{1-\omega}M[f]>\nu_\theta[f]M[f]>0.
\]
Thus the gain is bounded by $0<G_\theta[f]<C[f]$. Writing $p=\nu_\theta[f]M[f]/C[f]\in(0,1)$, monotonicity of $\sigma$ and the bound $|\tanh|\leq1$ give
\[
  G_\theta[f]\geq C[f]\sigma(\operatorname{logit}p-A_g)
  =\frac{\nu_\theta[f]M[f]}{p+(1-p)e^{A_g}}
  \geq e^{-A_g}r_{\min}\rho[f]T[f]^{1-\omega}M[f].
\]
The positive lower bound on $M[f]$ gives the required gain lower bound. At equilibrium, $g=0$ makes both perturbations vanish, so $\nu_\theta[M]=r_{\mathrm B}\rho[M]T[M]^{1-\omega}$. The identity $\sigma(\operatorname{logit}p)=p$ yields
\[
  G_\theta[M]=C[M]\sigma\left(\operatorname{logit}\frac{\nu_\theta[M]M}{C[M]}\right)=\nu_\theta[M]M.
\]

Under $f\mapsto af$, the temperature, encoding, normalized velocity, gate, and two network outputs are unchanged. The two terms inside the envelope scale by $a$, so
\[
  \nu_\theta[af]=a\nu_\theta[f],\quad C[af]=a^2C[f],\quad \frac{\nu_\theta[af]M[af]}{C[af]}=\frac{\nu_\theta[f]M[f]}{C[f]}.
\]
Substitution into~\eqref{eq:bounded-gain} gives $G_\theta[af]=a^2G_\theta[f]$. Finally, when $a_\theta[f]=b_\theta[f]=0$, the readouts give $\nu_\theta[f]=r_{\mathrm B}\rho[f]T[f]^{1-\omega}$ and $G_\theta[f]=\nu_\theta[f]M[f]$. Hence $q_\theta[f]=r_{\mathrm B}\rho[f]T[f]^{1-\omega}(M[f]-f)=Q_{\mathrm B}[f]$.
\end{proof}

We obtain $\mathcal{Q}_{\theta,\kappa}[f]$ by substituting the learned gain and frequency into~\eqref{eq:corrected-gain-loss}. The resulting operator inherits the structural and positivity guarantees of Section~\ref{sec:method} under the stated conditions.

\subsection{Discrete realization on the reference velocity grid}
\label{sec:discrete-realization}

The neural operator is discretized and trained on a fixed reference velocity grid in $\Omega_L=[-L,L]^3$. To retain the structural identities of Section~\ref{sec:method}, we use the same positive quadrature throughout the construction. Let $\{(\mathbf{v}_j,\omega_j)\}_{j=1}^{N_v}$ satisfy $\omega_j>0$, and write $\psi_j=\psi(\mathbf{v}_j)$. For a positive nodal distribution $f=(f_j)_{j=1}^{N_v}$, the density is $\rho_h[f]=\sum_j\omega_jf_j$, and the discrete equilibrium is determined by
\begin{equation}
  M_j[f]=e^{\alpha_h[f]^T\psi_j},\quad \sum_{j=1}^{N_v}\omega_j\psi_jM_j[f]=\sum_{j=1}^{N_v}\omega_j\psi_jf_j.
  \label{eq:discrete-equilibrium}
\end{equation}
This is the entropy-minimizing equilibrium used in discrete-velocity BGK models~\cite{discrete_bgk_model}. The same equilibrium is used in the encoding, gate, gain parameterization, and BGK reference. All velocity integrals and $L^2$ norms are evaluated using the same quadrature, while $L^\infty$ quantities are evaluated by maxima over the velocity nodes.

For the conservative projection, set $H_{f,h}=\sum_j\omega_jf_j\psi_j\psi_j^T$ and define
\begin{equation}
  (\mathcal{P}_{f,h}q)_j=q_j-\left(\psi_j^TH_{f,h}^{-1}\sum_{k=1}^{N_v}\omega_k\psi_kq_k\right)f_j.
  \label{eq:discrete-projection}
\end{equation}
Assume that the nodal invariant vectors $\{\psi_j\}_{j=1}^{N_v}$ span $\mathbb{R}^5$. Since $\omega_jf_j>0$, $H_{f,h}$ is positive definite, and the five discrete collision-invariant moments of the projected field vanish. Together with~\eqref{eq:discrete-equilibrium}, this ensures moment conservation for both the corrected proposal and the discrete BGK reference.

The discrete entropy pairing and relative entropy are
\begin{equation}
  D_{f,h}(q)=\sum_j\omega_jq_j\log\frac{f_j}{\rho_h[f]},\quad \mathscr{H}_h(f\mid M[f])=\sum_j\omega_jf_j\log\frac{f_j}{M_j[f]}.
  \label{eq:discrete-moments-entropy}
\end{equation}
Using these quantities in the entropy correction, the arguments of Section~\ref{sec:method} carry over with integrals replaced by weighted sums, giving the following bound with $T_h[f]$ evaluated using the same quadrature
\[
  D_{f,h}(\mathcal{Q}_{\theta,\kappa,h}[f])\leq-\kappa r_{\mathrm B}\rho_h[f]T_h[f]^{1-\omega}\mathscr{H}_h(f\mid M[f]).
\]
The equilibrium and density-scaling identities carry over in the same way, while the forward Euler positivity condition uses a maximum over the velocity nodes. These statements hold in exact arithmetic under exact moment matching.

\subsection{Training}
\label{sec:training}

We determine $\theta$ by fitting the corrected operator on the reference grid to Fourier spectral evaluations of the collision term. Let
\begin{equation}
  \mathcal{D}=\left\{\left(f^{(n)},\mathcal{Q}_{\mathrm{ref}}[f^{(n)},f^{(n)}]\right)\right\}_{n=1}^{N_{\mathrm{train}}},
  \label{eq:training_data}
\end{equation}
where $\mathcal{Q}_{\mathrm{ref}}$ is the spectral approximation from Section~\ref{sec::spectral} and $N_{\mathrm{train}}$ is the number of training samples. Appendix~\ref{app:dataset-generation} describes the sampling of the input distributions and the generation of their reference collision terms.

For a fixed $\kappa\in[0,1]$, we minimize the normalized error in the corrected collision operator,
\begin{equation}
  \mathcal{L}(\theta)=\frac{1}{N_{\mathrm{train}}}\sum_{n=1}^{N_{\mathrm{train}}}\frac{\|\mathcal{Q}_{\theta,\kappa,h}[f^{(n)}]-Q_{\mathrm{ref}}^{(n)}\|_{\mathrm{rms}}^2}{\|Q_{\mathrm{ref}}^{(n)}\|_{\mathrm{rms}}^2+\varepsilon_{\mathrm{loss}}^2},
  \label{eq:training-loss}
\end{equation}
where $Q_{\mathrm{ref}}^{(n)}$ denotes the reference values at the nodes, $\|q\|_{\mathrm{rms}}^2=N_v^{-1}\sum_j|q_j|^2$, and $\varepsilon_{\mathrm{loss}}>0$ regularizes the denominator near equilibrium. No auxiliary structural penalties are needed because the constraints are imposed by the construction. To start training from the BGK reference, we initialize the weights and biases of both affine readouts to zero, as justified by Proposition~\ref{prop:neural-realization}. The numerical settings are given in Section~\ref{sec:numerical-setup}.

\subsection{Transfer to general velocity domains}
\label{sec:velocity-domain-transfer}

The operator is trained on the fixed reference box $\Omega_L=[-L,L]^3$, whereas a kinetic solver may use $\Omega_{L_v}=\mathbf{v}_0+[-L_v,L_v]^3$ with $L_v>0$. To reuse the trained model, we map the physical velocity grid to the reference box and exploit the translation invariance and velocity-scaling covariance of the collision operator. Let
\begin{equation}
  s=\frac{L}{L_v},\quad \widehat{\mathbf{v}}=s(\mathbf{v}-\mathbf{v}_0),\quad \widehat f(\widehat{\mathbf{v}})=f(\mathbf{v}_0+\widehat{\mathbf{v}}/s).
  \label{eq:velocity-domain-pullback}
\end{equation}
For corresponding uniform Cartesian grids with the same number of nodes, $\widehat f_j=f_j$, so the network input is unchanged.

For the VHS kernel $\mathcal{B}(\mathbf{g},\boldsymbol{\sigma})=C_\omega|\mathbf{g}|^{2(1-\omega)}$ in~\eqref{eq:vhs_kernel}, a change of variables in the collision integral using $\widehat{\mathbf{v}}-\widehat{\mathbf{v}}_*=s(\mathbf{v}-\mathbf{v}_*)$ gives
\begin{equation}
  \mathcal{Q}[f,f](\mathbf{v})=s^{-(5-2\omega)}\mathcal{Q}[\widehat f,\widehat f](s(\mathbf{v}-\mathbf{v}_0)).
  \label{eq:collision-velocity-scaling}
\end{equation}
Suppressing the fixed parameters $\theta,\kappa$ and the discretization index $h$, write $\mathcal{Q}^{L}$ for the trained reference-domain operator. We define
\begin{equation}
  \bigl(\mathcal{Q}^{L_v}[f]\bigr)_j=s^{-(5-2\omega)}\bigl(\mathcal{Q}^{L}[\widehat f]\bigr)_j,\quad \widehat{\mathbf{v}}_j=s(\mathbf{v}_j-\mathbf{v}_0).
  \label{eq:learned-velocity-interface}
\end{equation}
The solver therefore reuses the trained model without retraining or interpolation and rescales only its output. For the constant kernel with $\omega=1$, the factor is $s^{-3}$.

The transfer is compatible with the structural construction. The moments satisfy $\rho[\widehat f]=s^3\rho[f]$, $\mathbf{u}[\widehat f]=s(\mathbf{u}[f]-\mathbf{v}_0)$, and $T[\widehat f]=s^2T[f]$, while the normalized invariant vectors agree in the centered coordinates $\widehat{\mathbf{v}}_j/L=(\mathbf{v}_j-\mathbf{v}_0)/L_v$. Centering preserves the span of the collision invariants. Together with the quadrature scaling $\widehat\omega_j=s^3\omega_j$, these relations preserve moment matching, conservation, and equilibrium balance. Mass conservation cancels the constant shift in the normalized logarithm, so the physical entropy pairing is $s^{-(8-2\omega)}$ times the reference pairing and retains its sign.

\section{Numerical experiments}
\label{sec:experiments}

\subsection{Numerical setup}
\label{sec:numerical-setup}

We consider the normalized constant collision kernel $\mathcal{B}=\frac{1}{4\pi}$. Unless otherwise stated, we use a uniform $64^3$ Fourier grid on the periodic cell $[-L,L)^3$ with $L=10$ and relative-velocity cutoff $R=9$, and perform spectral evaluations in FP64. In the numerical experiments, $p$ denotes the degree of exactness of the Lebedev angular quadrature. Reference solutions are analytic when available and otherwise computed by the spectral method with $p=23$. Relative $L^1$ and $L^2$ errors are evaluated using the uniform velocity quadrature.

The learned operator is implemented in PyTorch~\cite{pytorch} and trained on this reference grid. We set $A_g=A_r=10$, $r_{\min}=0.02$, $r_{\max}=8$, $r_{\mathrm B}=0.5$, and $\tau=0.1$. All reported experiments use $\kappa=0$, which enforces the non-strict entropy inequality. We train with $\varepsilon_{\mathrm{loss}}=10^{-4}$ in~\eqref{eq:training-loss}, using AdamW~\cite{adamw} with a cosine learning-rate schedule~\cite{cos_learning} following an initial warmup period. The main training and computationally intensive experiments were performed on an NVIDIA A40 GPU. All reported runtime measurements were obtained on the same device.

\subsection{Accuracy and computational cost}
\label{sec:operator-accuracy}

We first use a fixed nonequilibrium two-stream distribution to illustrate the accuracy and computational cost of a single collision evaluation. We take
\[
f(\mathbf{v})=\frac{1}{2(2\pi T_0)^{3/2}}\left[\exp\left(-\frac{|\mathbf{v}-\mathbf{u}_0|^2}{2T_0}\right)+\exp\left(-\frac{|\mathbf{v}+\mathbf{u}_0|^2}{2T_0}\right)\right],
\]
where $\mathbf{u}_0=(2.5,0,0)$ and $T_0=0.6$, and normalize its nodal values to unit mass on the $64^3$ grid. For the $32^3$ comparisons, the input and spectral reference are restricted to the coincident velocity nodes. Per-sample times are GPU-synchronized wall-clock times for batches of four samples, divided by four and averaged over 20 repetitions after five warmup calls, excluding setup and data transfers.

Tables~\ref{tab:collision-cost} and~\ref{tab:collision-cost-n64} show that, for this input, the FP32 learned operator is both faster and more accurate than the spectral method with $p=7$ and $9$ on both grids. Relative to the $p=9$ spectral evaluation, the learned operator reduces the relative $L^2$ error from about $6.5\%$ to below $5\%$, while giving speedups of $12.80\times$ and $29.86\times$ on the $32^3$ and $64^3$ grids, respectively. Increasing $p$ to $11$ yields substantially smaller errors at higher computational cost. The FP64 and FP32 errors of the learned operator are close, with differences below the precision reported in the tables. We therefore use FP32 for evaluations of the learned operator in the following experiments unless otherwise stated.

\begin{table}[pos=htbp]
\caption{Two-stream collision evaluation on the $32^3$ grid. Spectral 7, 9, and 11 refer to $p=7,9,11$, respectively. Times are mean $\pm$ standard deviation. Speedup is the spectral mean time divided by the FP32 mean time of the learned operator on this grid.}
\label{tab:collision-cost}
\centering
\setlength{\tabcolsep}{3pt}
\begin{tabular}{lccccc}
\hline
Method & Precision & Rel. $L^1$ [\%] & Rel. $L^2$ [\%] & Time [ms] & Speedup \\
\hline
\textbf{Learned operator} & \textbf{FP32} & $\mathbf{5.306}$ & $\mathbf{4.827}$ & $\mathbf{1.541\pm0.007}$ & --- \\
Learned operator & FP64 & $5.306$ & $4.827$ & $4.873\pm0.025$ & --- \\
\hline
Spectral 7 & FP64 & $150.5$ & $49.36$ & $15.789\pm0.142$ & $10.25\times$ \\
Spectral 9 & FP64 & $19.32$ & $6.494$ & $19.715\pm0.025$ & $12.80\times$ \\
Spectral 11 & FP64 & $0.8221$ & $0.2518$ & $26.839\pm0.127$ & $17.42\times$ \\
\hline
\end{tabular}
\end{table}

\begin{table}[pos=htbp]
\caption{Two-stream collision evaluation on the $64^3$ grid, using the conventions of Table~\ref{tab:collision-cost}. Speedups use the FP32 timing of the learned operator on the $64^3$ grid.}
\label{tab:collision-cost-n64}
\centering
\setlength{\tabcolsep}{3pt}
\begin{tabular}{lccccc}
\hline
Method & Precision & Rel. $L^1$ [\%] & Rel. $L^2$ [\%] & Time [ms] & Speedup \\
\hline
\textbf{Learned operator} & \textbf{FP32} & $\mathbf{5.391}$ & $\mathbf{4.918}$ & $\mathbf{5.256\pm0.004}$ & --- \\
Learned operator & FP64 & $5.391$ & $4.918$ & $27.205\pm0.008$ & --- \\
\hline
Spectral 7 & FP64 & $150.7$ & $49.35$ & $109.211\pm0.008$ & $20.78\times$ \\
Spectral 9 & FP64 & $19.36$ & $6.495$ & $156.953\pm0.012$ & $29.86\times$ \\
Spectral 11 & FP64 & $0.7912$ & $0.2387$ & $207.239\pm0.013$ & $39.43\times$ \\
\hline
\end{tabular}
\end{table}

\begin{figure}[pos=htbp]
\centering
\includegraphics[width=\textwidth]{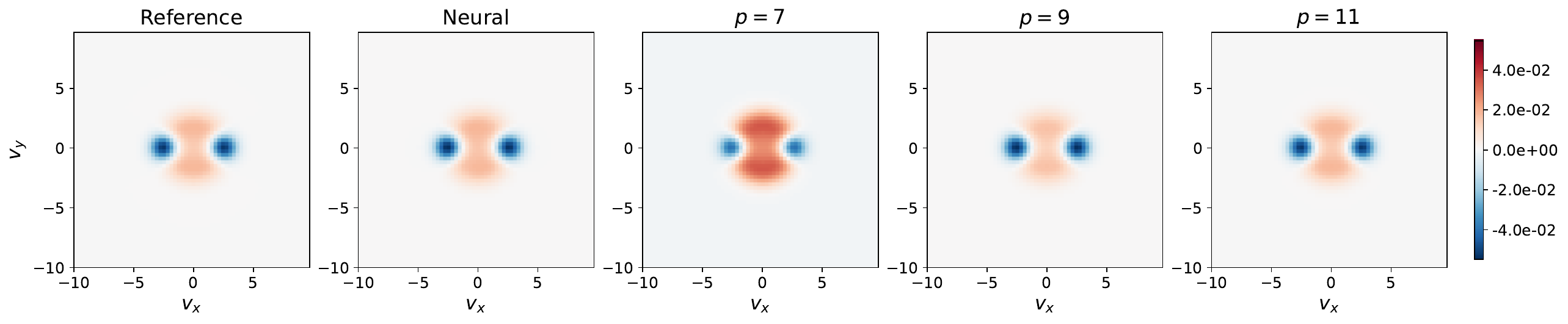}
\caption{Collision terms for the two-stream distribution, integrated along $v_z$. From left to right: the spectral reference, the prediction of the learned operator, and spectral evaluations with $p=7,9,11$. All panels use the same color scale.}
\label{fig:two-stream-collision}
\end{figure}

Figure~\ref{fig:two-stream-collision} compares the $v_z$-integrated collision fields $\int_{-L}^{L}q(v_x,v_y,v_z)\,\mathrm{d}v_z$. The spectral result at $p=7$ shows visible distortion, whereas the prediction of the learned operator and the spectral results at $p=9$ and $11$ are visually close to the reference.

To assess accuracy beyond a single input, we also evaluate the learned operator and the BGK approximation $\frac{1}{2}\rho[f](M[f]-f)$ on the 1000 held-out test samples from the dataset described in Appendix~\ref{app:dataset-generation}, using FP32 and the same discrete moment-matched equilibrium. Table~\ref{tab:test-set-collision-error} reports relative $L^2$ errors in the collision term against the spectral reference. The learned operator reduces both the mean error and its variation across samples compared with BGK.

\begin{table}[pos=htbp]
\caption{Collision-term relative $L^2$ errors on 1000 test samples. Both evaluations use FP32. The sample standard deviation uses denominator $N_{\mathrm{test}}-1$.}
\label{tab:test-set-collision-error}
\centering
\begin{tabular}{lcc}
\hline
Method & Mean [\%] & Std. dev. [\%] \\
\hline
Learned operator & $4.627$ & $4.165$ \\
BGK & $28.326$ & $13.597$ \\
\hline
\end{tabular}
\end{table}

The following experiments examine the resulting collision dynamics and kinetic solutions.

\subsection[Homogeneous Problems]{Homogeneous Problems}
\label{sec:homogeneous-relaxation}

We next consider the spatially homogeneous Boltzmann equation
\begin{equation}
\partial_t f=\mathcal{Q}[f,f],\quad f(0,\mathbf{v})=f_0(\mathbf{v}).
\label{eq:homogeneous-test}
\end{equation}
All trajectories are advanced with the second-order SSP Runge--Kutta method~\cite{SSPRK} using $\Delta t=0.1$.

\medskip
\noindent\textbf{BKW.}
Our first test uses the analytic Bobylev--Krook--Wu (BKW) solution~\cite{BKW}. For the constant kernel $\mathcal{B}=1/(4\pi)$,
\begin{equation}
f_{\mathrm{BKW}}(t,\mathbf{v})=\frac{e^{-|\mathbf{v}|^2/(2K(t))}}{(2\pi K(t))^{3/2}}\left[\frac{5K(t)-3}{2K(t)}+\frac{1-K(t)}{2K(t)^2}|\mathbf{v}|^2\right],
\label{eq:bkw-solution}
\end{equation}
where $K(t)=1-(1-K_0)e^{-t/6}$. For $3/5\leq K_0\leq1$, this is a nonnegative unit-mass solution on $\mathbb{R}^3$ with zero mean velocity and unit temperature that relaxes to $\mathcal{M}_{1,\mathbf{0},1}$. We set $K_0=0.6$ and evolve the sampled initial data to $t=2$. For comparison, we use the BGK model $\partial_t f=\tfrac12\rho[f](M[f]-f)$ with the discrete moment-matched equilibrium $M[f]$.

\begin{figure}[pos=htbp]
\centering
\subfloat[$t=1$]{\includegraphics[width=0.32\textwidth]{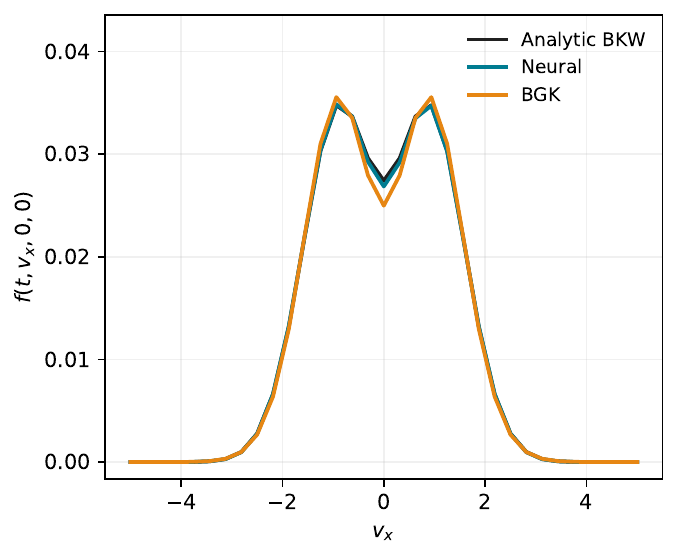}}\hfill
\subfloat[$t=2$]{\includegraphics[width=0.32\textwidth]{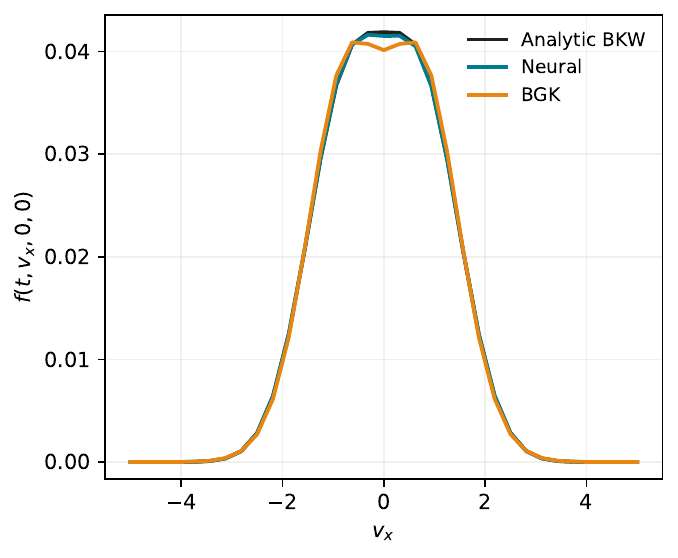}}\hfill
\subfloat[Relative error]{\includegraphics[width=0.32\textwidth]{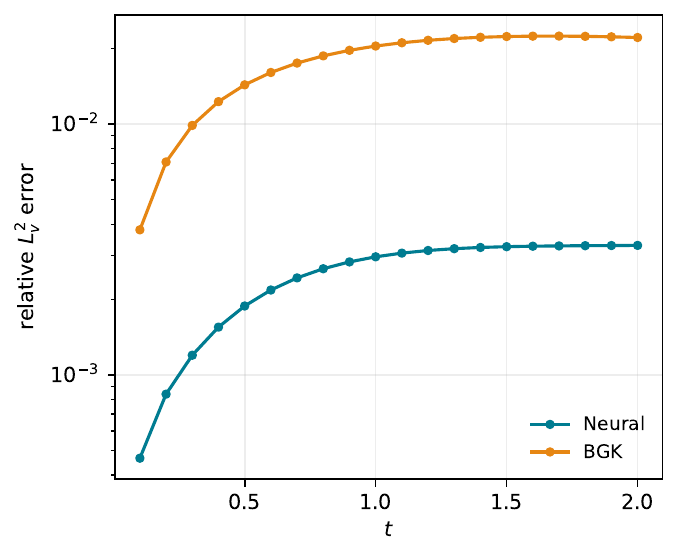}}
\caption{BKW relaxation with $K_0=0.6$. Left and center: profiles $f(t,v_x,0,0)$ at $t=1$ and $2$ from the analytic solution, the learned operator, and BGK. Right: relative $L^2$ errors against the analytic solution.}
\label{fig:bkw-relaxation}
\end{figure}

Figure~\ref{fig:bkw-relaxation} shows that the solution obtained with the learned operator closely follows the analytic BKW evolution. Although the BGK frequency $\rho[f]/2$ matches the deviatoric-pressure relaxation rate of the constant kernel, its single-rate Maxwellian relaxation underestimates the profile near $v_x=0$ at both times. It captures this structure more accurately and has lower relative $L^2$ error at every recorded time in $0<t\leq2$.

\medskip
\noindent\textbf{Three-stream.}
We next consider the asymmetric three-stream distribution
\begin{equation}
f_0(\mathbf{v})=\sum_{i=1}^{3}\frac{w_i}{(2\pi T_0)^{3/2}}\exp\left(-\frac{|\mathbf{v}-\mathbf{u}_i|^2}{2T_0}\right),
\label{eq:three-stream-initial}
\end{equation}
with $T_0=0.6$, weights $(w_1,w_2,w_3)=(0.4,0.35,0.25)$, and velocities
\[
\mathbf{u}_1=(2.5,0,0),\quad \mathbf{u}_2=(-1.6,1.4,0),\quad \mathbf{u}_3=(-1.76,-1.96,0).
\]
Each Maxwellian component is discretely normalized to unit mass before forming the mixture.

We evolve the system using the learned operator and the spectral method with $p=9$ to $t=5$ and compare both against an independently evolved spectral reference.

\begin{figure}[pos=htbp]
\centering
\subfloat[Collision field]{\includegraphics[width=0.49\textwidth]{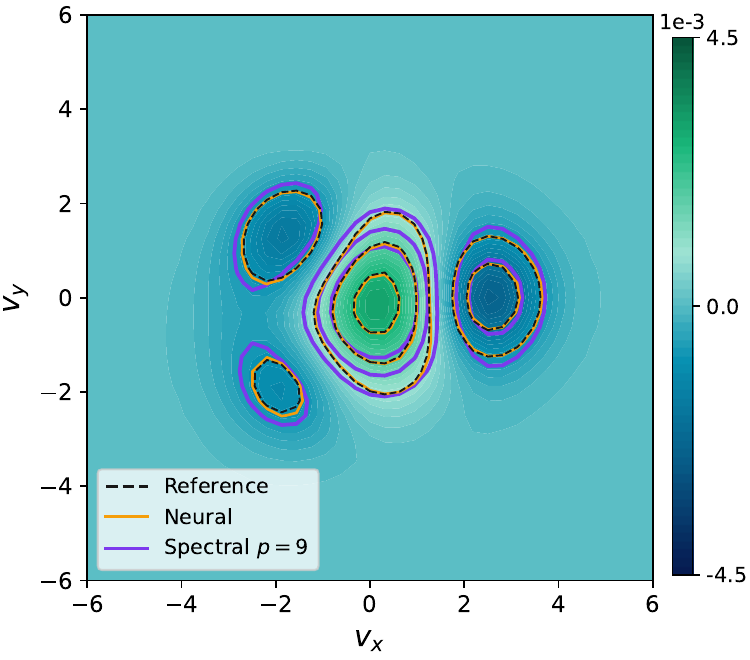}}\hfill
\subfloat[Relative error]{\includegraphics[width=0.49\textwidth]{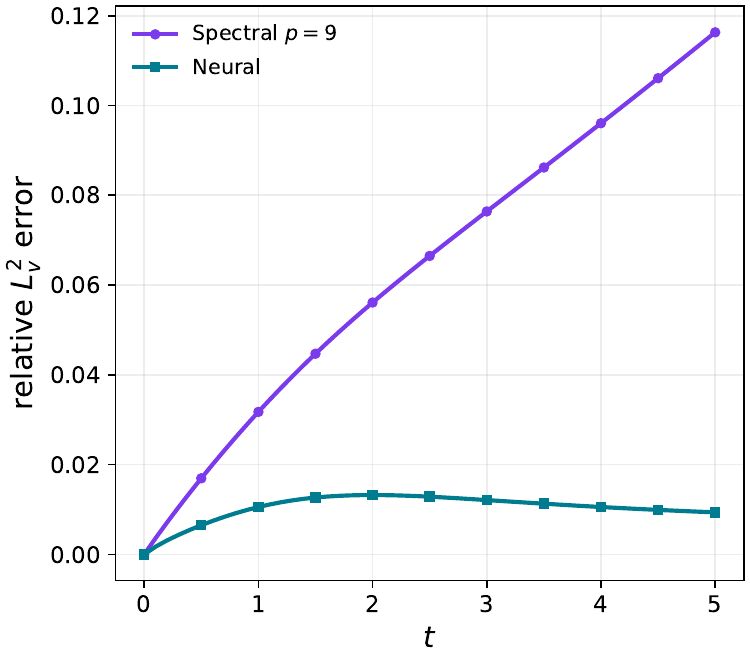}}
\caption{Three-stream relaxation. Left: $v_z$-integrated collision fields at $t=5$, evaluated by each method on its own evolved state. Colors show the collision field obtained with the learned operator, and contours compare the three methods. Right: relative $L^2$ errors of the evolved distributions against the spectral reference.}
\label{fig:three-stream-relaxation}
\end{figure}

Figure~\ref{fig:three-stream-relaxation} shows that the collision field obtained with the learned operator at $t=5$ remains close to the reference, while the $p=9$ spectral result shows noticeable deviations in the contour geometry. The resulting trajectory also maintains a lower relative $L^2$ error in $f$ over the simulated interval.

\begin{figure}[pos=htbp]
\centering
\subfloat[Conservation]{\includegraphics[width=0.32\textwidth]{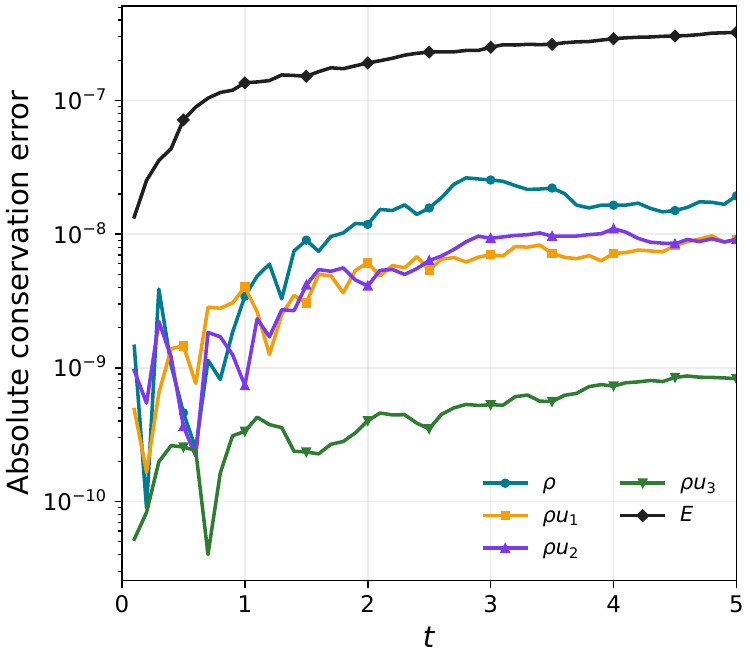}}\hfill
\subfloat[Minimum value]{\includegraphics[width=0.32\textwidth]{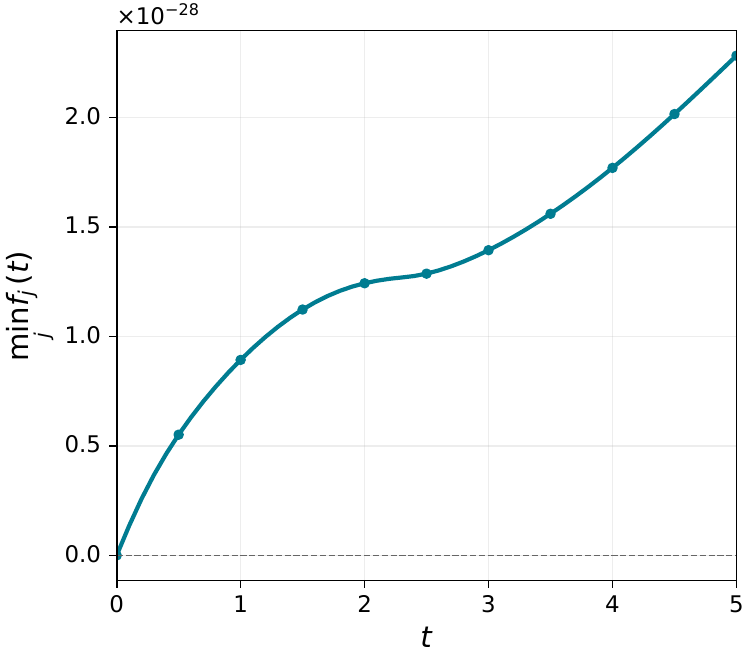}}\hfill
\subfloat[Relative entropy]{\includegraphics[width=0.32\textwidth]{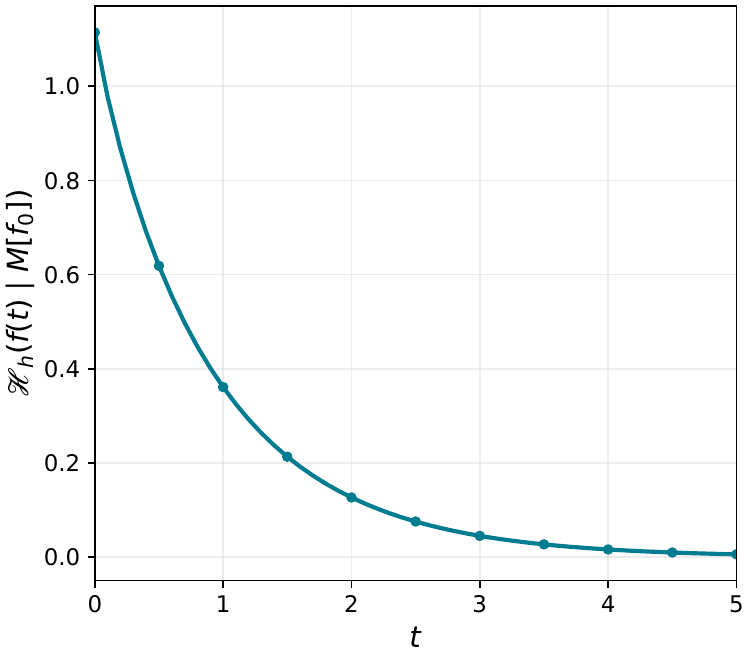}}
\caption{Three-stream structural diagnostics. Left: absolute conservation errors. Center: minimum nodal value. Right: relative entropy with respect to $M[f_0]$.}
\label{fig:three-stream-structure}
\end{figure}

To verify conservation, positivity, and entropy decay along the learned trajectory, Figure~\ref{fig:three-stream-structure} reports the corresponding structural diagnostics. The largest absolute drifts in mass, individual momentum components, and energy are $2.64\times10^{-8}$, $1.10\times10^{-8}$, and $3.24\times10^{-7}$, respectively, consistent with roundoff from FP32 collision evaluations. The distribution remains nonnegative at all time levels and Runge--Kutta stages. The discrete relative entropy decreases monotonically at every time step, from $1.114$ initially to $5.95\times10^{-3}$ at $t=5$.

\medskip
\noindent\textbf{Discontinuous initial data.}
To assess generalization beyond the smooth training distribution, we consider a nonsmooth input obtained by introducing a jump across $v_y=0$ into the two-stream distribution
\[
g(\mathbf{v})=\frac{1}{2}\bigl(\mathcal{M}_{1,\mathbf{u}_-,T_0}(\mathbf{v})+\mathcal{M}_{1,\mathbf{u}_+,T_0}(\mathbf{v})\bigr),
\]
with $\mathbf{u}_\pm=(\pm2.5,0,0)$ and $T_0=0.6$. We set
\begin{equation}
f_0(\mathbf{v})=C g(\mathbf{v})\left[(1-\varepsilon)\mathbf{1}_{\{v_y<0\}}+(1+\varepsilon)\mathbf{1}_{\{v_y\geq0\}}\right],\quad \varepsilon=0.25,
\label{eq:discontinuous-initial}
\end{equation}
where $C$ normalizes the sampled initial distribution to unit mass. We apply the same trained operator without retraining and evolve the system to $t=5$, comparing against the spectral method with $p=9$ and the independently evolved spectral reference.

\begin{figure}[pos=htbp]
\centering
\subfloat[Initial distribution]{\includegraphics[width=0.32\textwidth]{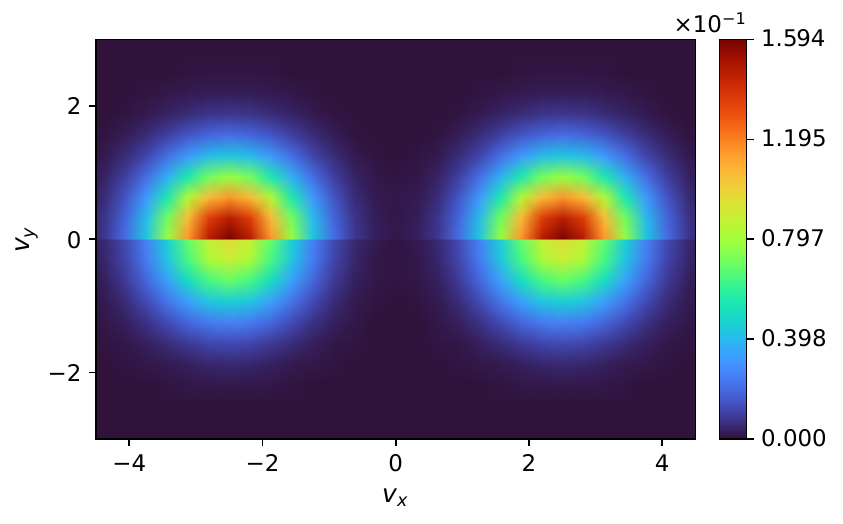}}\hfill
\subfloat[Reference collision]{\includegraphics[width=0.32\textwidth]{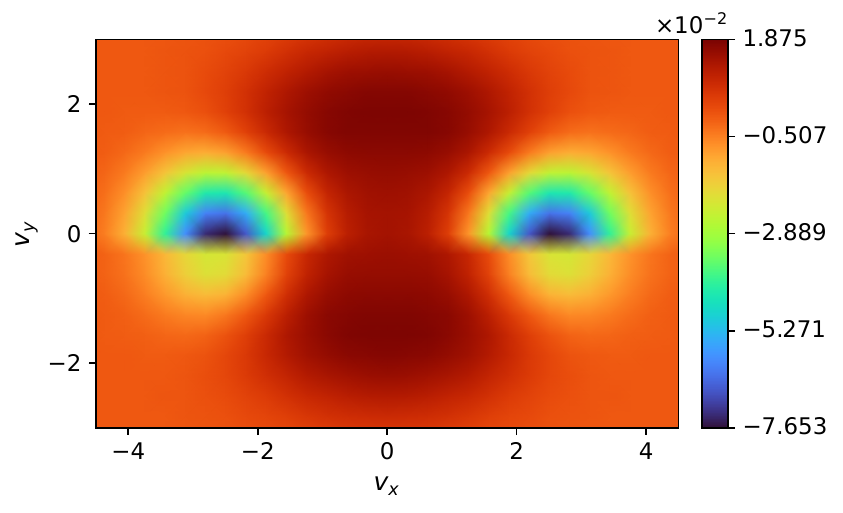}}\hfill
\subfloat[Learned collision]{\includegraphics[width=0.32\textwidth]{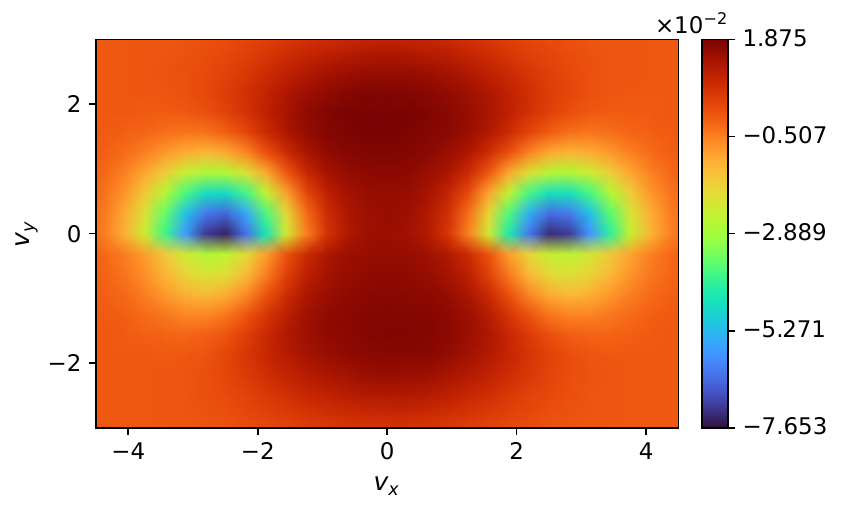}}
\caption{Discontinuous two-stream initial data with $\varepsilon=0.25$ and the corresponding reference and learned collision terms. All fields are integrated over $v_z$, and the two collision panels use the same color scale.}
\label{fig:discontinuous-initial}
\end{figure}

\begin{figure}[pos=htbp]
\centering
\subfloat[Collision field]{\includegraphics[width=0.32\textwidth]{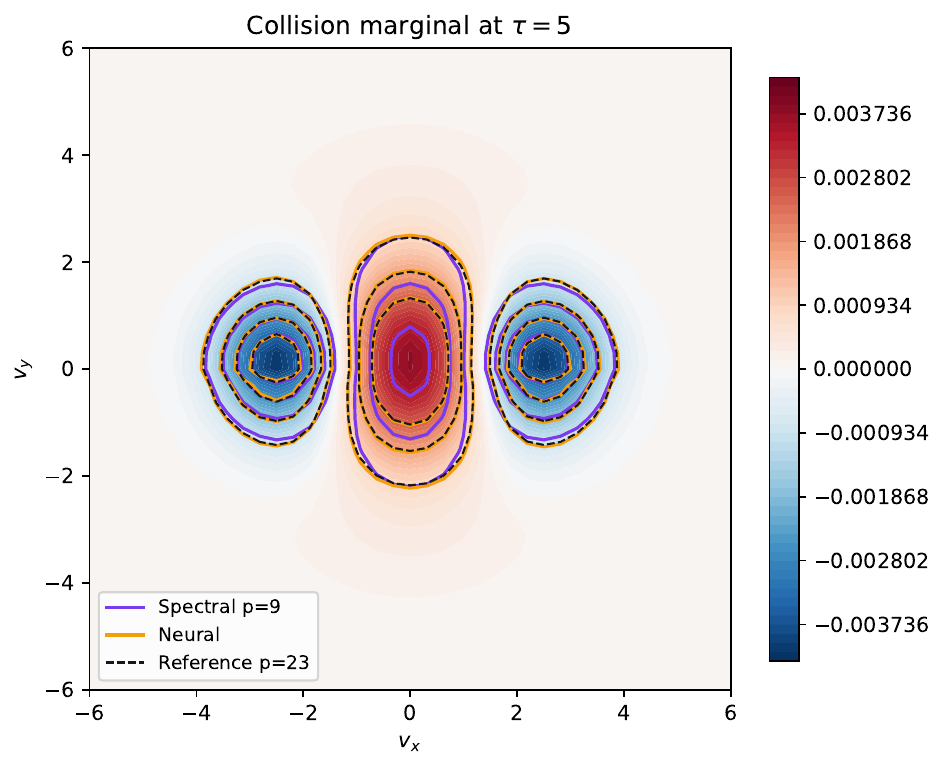}}\hfill
\subfloat[$L^1$ error]{\includegraphics[width=0.32\textwidth]{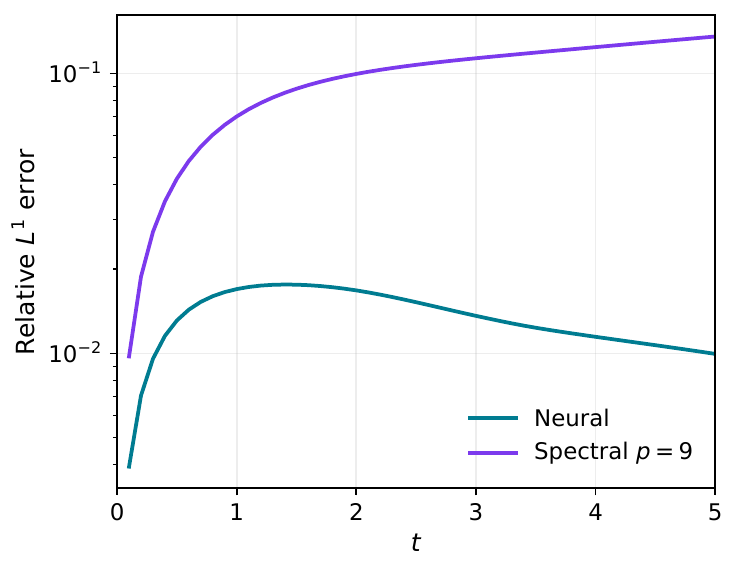}}\hfill
\subfloat[$L^2$ error]{\includegraphics[width=0.32\textwidth]{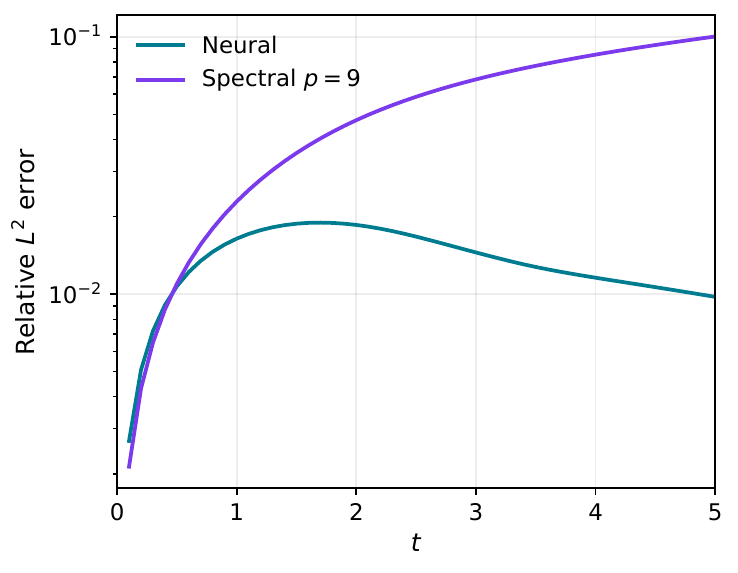}}
\caption{Relaxation from discontinuous initial data. Left: $v_z$-integrated collision fields at $t=5$, evaluated by each method on its own evolved state, with the reference shown in color. Center and right: relative errors in $f$ against the spectral reference.}
\label{fig:discontinuous-relaxation}
\end{figure}

Figure~\ref{fig:discontinuous-initial} shows that the learned operator captures the asymmetric collision structure induced by the jump, including the two localized loss regions. Figure~\ref{fig:discontinuous-relaxation} compares the subsequent evolution. At $t=5$, the collision contours obtained with the learned operator remain close to the reference, while the $p=9$ spectral contours exhibit larger deviations. The relative $L^1$ and $L^2$ errors of the solution obtained with the learned operator decrease after an initial rise and reach approximately $1\%$ at $t=5$, compared with approximately $13.5\%$ and $10.0\%$ for the spectral solution. This test demonstrates that the operator trained on smooth distributions also accurately resolves relaxation from the discontinuous input considered here.

\subsection[Spatially Inhomogeneous Problems]{Spatially Inhomogeneous Problems}
\label{sec:steady-problems}
\label{sec:inhomogeneous-dynamics}

We next consider spatially inhomogeneous kinetic problems. For Couette flow, we solve the stationary equation by direct iteration~\cite{linearized_boltzmann_collision_operator}. For the normal shock and lid-driven cavity, we advance the time-dependent equation
\begin{equation}
\partial_t f+\sum_{i=1}^{d_x}v_i\partial_{x_i}f=\frac{1}{\mathrm{Kn}}\mathcal{Q}[f,f],\quad d_x\in\{1,2\},
\label{eq:inhomogeneous-tests}
\end{equation}
and compare the macroscopic fields at the final time.

\medskip
\noindent\textbf{Couette flow.}

We consider the steady Boltzmann equation in one spatial dimension,
\begin{equation}
v_1\partial_x f=\frac{1}{\mathrm{Kn}}\mathcal{Q}[f,f],\quad x\in\left(-\frac{1}{2},\frac{1}{2}\right).
\label{eq:steady-couette}
\end{equation}
Both walls have temperature $T_{\mathrm w}=1$. Their tangential velocities are $\mathbf{u}_{\mathrm L}=(0,-0.5,0)$ and $\mathbf{u}_{\mathrm R}=(0,0.5,0)$. We impose diffuse Maxwellian boundary conditions on the incoming velocities,
\begin{equation}
\begin{cases}
f\left(-\frac{1}{2},\mathbf{v}\right)=\rho_{\mathrm L}\mathcal{M}_{1,\mathbf{u}_{\mathrm L},T_{\mathrm w}}(\mathbf{v}), & v_1>0,\\
f\left(\frac{1}{2},\mathbf{v}\right)=\rho_{\mathrm R}\mathcal{M}_{1,\mathbf{u}_{\mathrm R},T_{\mathrm w}}(\mathbf{v}), & v_1<0.
\end{cases}
\label{eq:couette-wall}
\end{equation}
The wall densities $\rho_{\mathrm L}$ and $\rho_{\mathrm R}$ are determined by zero normal mass flux at each wall. We also prescribe the total mass,
\begin{equation}
\int_{-\frac{1}{2}}^{\frac{1}{2}}\int_{\Omega_L}f(x,\mathbf{v})\,\mathrm{d}\mathbf{v}\,\mathrm{d}x=1.
\label{eq:couette-total-mass}
\end{equation}

We consider $\mathrm{Kn}=0.1$, $1$, and $10$, using 32 spatial cells and a $64^3$ velocity grid. Transport is discretized by first-order upwinding, and the steady system is solved by a Newton-type iteration following~\cite{linearized_boltzmann_collision_operator}, with collision evaluations supplied by the learned operator. The wall densities are updated from the outgoing fluxes, and each outer iterate is normalized to preserve the prescribed total mass. We include a BGK baseline with the same boundary conditions and total mass.

\begin{figure}[pos=htbp]
\centering
\subfloat[Velocity]{\includegraphics[width=0.49\textwidth]{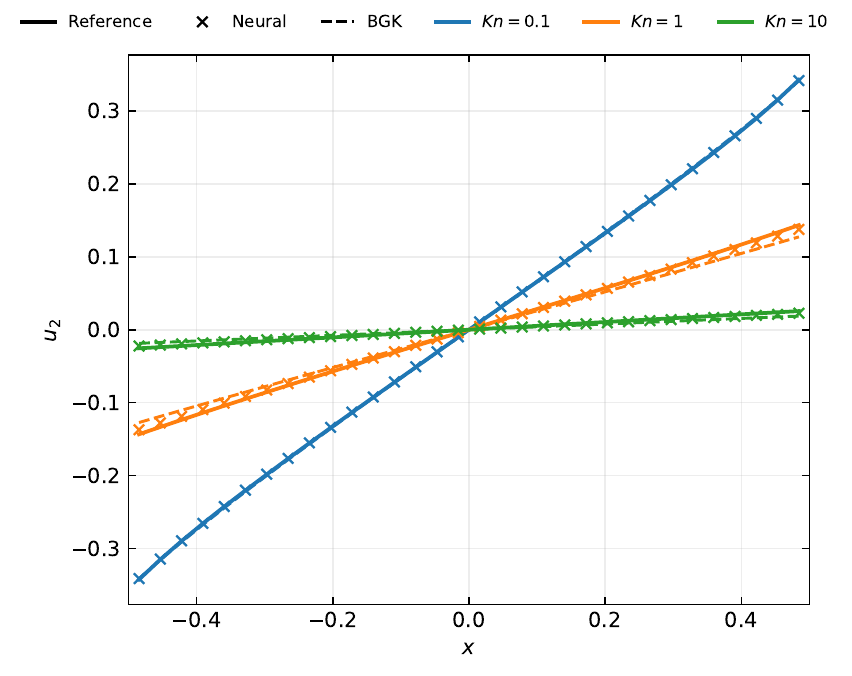}}\hfill
\subfloat[Temperature]{\includegraphics[width=0.49\textwidth]{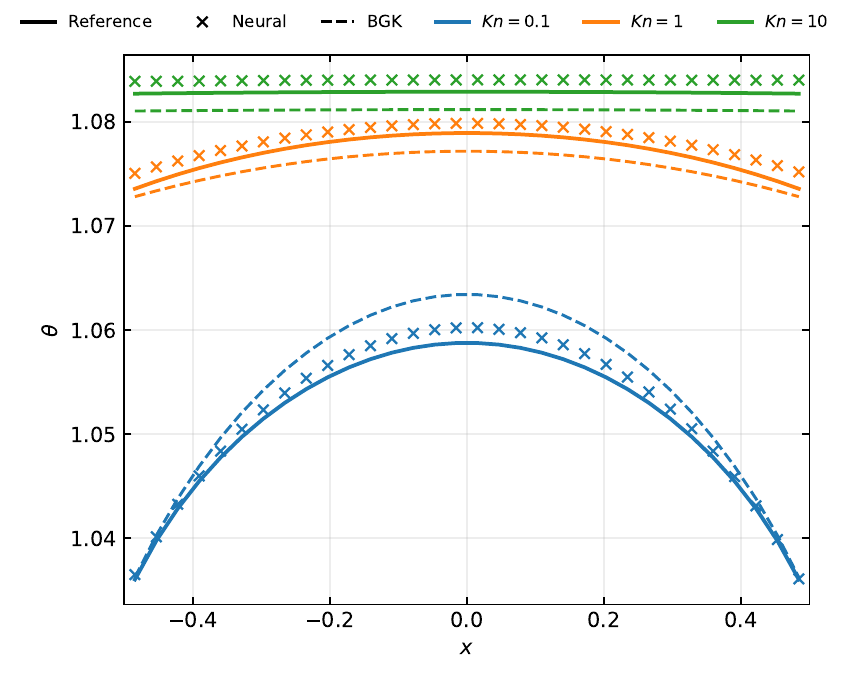}}
\caption{Steady Couette flow for $\mathrm{Kn}=0.1$, $1$, and $10$. Left: tangential velocity $u_2$. Right: temperature $\theta$. Solid lines, crosses, and dashed lines denote the spectral reference, solutions obtained with the learned operator, and steady BGK solutions, respectively.}
\label{fig:steady-couette}
\end{figure}

Figure~\ref{fig:steady-couette} shows that the solutions obtained with the learned operator agree well with the spectral reference across the three Knudsen numbers. Agreement improves over the BGK baseline, particularly for the temperature at $\mathrm{Kn}=0.1$ and the tangential velocity at $\mathrm{Kn}=1$.

For each spectral reference, we restart the Newton iteration from the converged solution obtained with the learned operator and switch to spectral collision evaluations.

\medskip
\noindent\textbf{Normal shock.}
We consider the normal-shock configuration in~\cite{unsteady_boltzmann_examples} on $x\in(-35,35)$, with initial data
\begin{equation}
f_0(x,\mathbf{v})=
\begin{cases}
\mathcal{M}_{\rho_-,\mathbf{u}_-,T_-}(\mathbf{v}), & x<0,\\
\mathcal{M}_{\rho_+,\mathbf{u}_+,T_+}(\mathbf{v}), & x\geq0,
\end{cases}
\quad \mathbf{u}_\pm=(u_\pm,0,0).
\label{eq:shock-initial}
\end{equation}
The upstream and downstream states satisfy the Rankine--Hugoniot relations
\begin{equation}
\begin{aligned}
\frac{\rho_+}{\rho_-}&=\frac{(\gamma+1)\mathrm{Ma}^2}{(\gamma-1)\mathrm{Ma}^2+2},\\
\frac{u_+}{u_-}&=\frac{(\gamma-1)\mathrm{Ma}^2+2}{(\gamma+1)\mathrm{Ma}^2},\\
\frac{T_+}{T_-}&=\frac{\bigl[(\gamma-1)\mathrm{Ma}^2+2\bigr]\bigl[2\gamma\mathrm{Ma}^2-\gamma+1\bigr]}{(\gamma+1)^2\mathrm{Ma}^2}.
\end{aligned}
\label{eq:shock-jump-density}
\end{equation}
We take $\mathrm{Ma}=2$, $\gamma=\frac{5}{3}$, $\rho_-=T_-=1$, and $u_-=\mathrm{Ma}\sqrt{\gamma T_-}$. The incoming distributions are fixed at the corresponding far-field Maxwellians,
\[
f(t,-35,\mathbf{v})=\mathcal{M}_{\rho_-,\mathbf{u}_-,T_-}(\mathbf{v})\quad(v_1>0),\quad
f(t,35,\mathbf{v})=\mathcal{M}_{\rho_+,\mathbf{u}_+,T_+}(\mathbf{v})\quad(v_1<0).
\]

We use $\mathrm{Kn}=0.1$ and advance to $t=50$ on 100 spatial cells with edges $x_j=35\,\sinh((j-50)/10)/\sinh(5)$, $j=0,\ldots,100$, concentrating the grid near the shock. The computation with the learned operator uses a $32^3$ velocity grid and conservative semi-Lagrangian transport with CWENO3 reconstruction~\cite{CWENO} and a positivity limiter. Strang splitting combines two half transport steps with a full SSPRK2 collision step, using $\Delta t=0.035$.

For comparison, DSMC~\cite{DSMC} uses the same spatial grid and $\Delta t\approx0.00236232$. The particle weight corresponds to 4000 upstream particles per uniform base cell of width $0.7$. The moments from 100 independent realizations are averaged at $t=50$ to obtain the macroscopic quantities.

\begin{figure}[pos=htbp]
\centering
\subfloat[Density]{\includegraphics[width=0.32\textwidth]{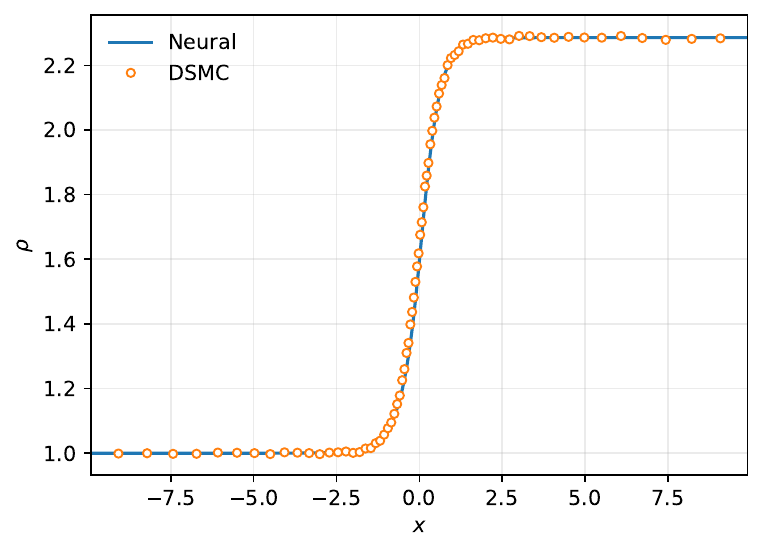}}\hfill
\subfloat[Velocity]{\includegraphics[width=0.32\textwidth]{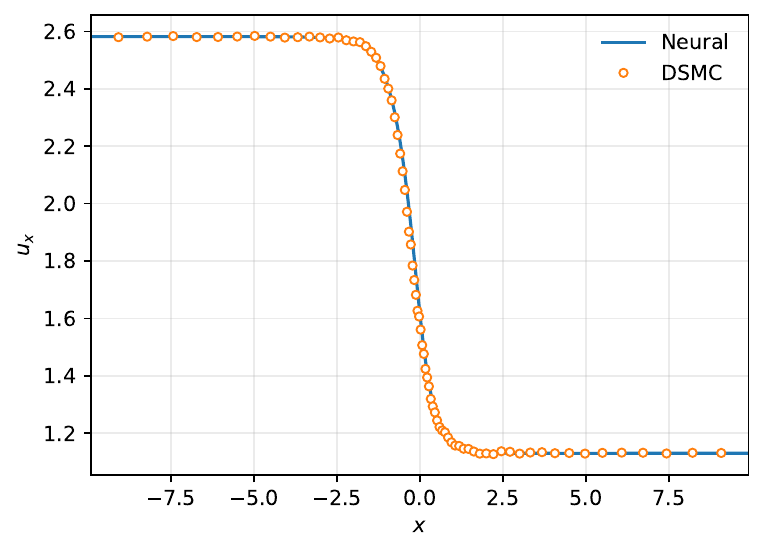}}\hfill
\subfloat[Temperature]{\includegraphics[width=0.32\textwidth]{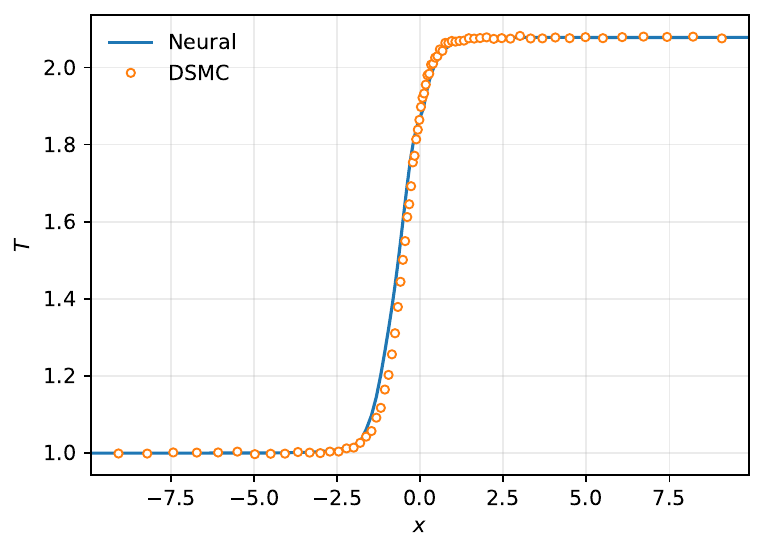}}
\caption{Normal shock at $\mathrm{Ma}=2$, $\mathrm{Kn}=0.1$, and $t=50$, displayed on $x\in(-10,10)$ to resolve the shock thickness. Solid lines denote the solution obtained with the learned operator, and circles denote the DSMC ensemble average.}
\label{fig:normal-shock}
\end{figure}

In Figure~\ref{fig:normal-shock}, the density, velocity, and temperature profiles agree closely with the DSMC ensemble averages. The solution obtained with the learned operator captures the shock location and thickness, with small differences within the shock layer. The deterministic solver produces smooth profiles in a single run, while mild sampling fluctuations remain after averaging 100 DSMC realizations. This test complements the single-evaluation benchmarks in Section~\ref{sec:operator-accuracy} by assessing the learned operator under repeated use in a time-dependent kinetic simulation.

\medskip
\noindent\textbf{Lid-driven cavity.}
On $(x,y)\in(0,1)^2$, the gas is initially at rest,
\begin{equation}
f_0(x,y,\mathbf{v})=\mathcal{M}_{1,\mathbf{0},1}(\mathbf{v}).
\label{eq:cavity-initial}
\end{equation}
All walls have temperature $T_{\mathrm w}=1$. The top lid has velocity $\mathbf{u}_{\mathrm w}=(1,0,0)$, while the other walls are stationary. Diffuse reflection is imposed on incoming velocities,
\begin{equation}
f(t,\mathbf{x},\mathbf{v})=\rho_{\mathrm w}(t,\mathbf{x})\mathcal{M}_{1,\mathbf{u}_{\mathrm w},T_{\mathrm w}}(\mathbf{v}),\quad \mathbf{v}\cdot\mathbf{n}<0,
\label{eq:cavity-wall}
\end{equation}
where $\mathbf{n}=(n_x,n_y,0)$ is the outward unit normal. Zero normal mass flux determines the wall density,
\[
\rho_{\mathrm w}=-\frac{\displaystyle\int_{\mathbf{v}\cdot\mathbf{n}>0}(\mathbf{v}\cdot\mathbf{n})f\,\mathrm{d}\mathbf{v}}{\displaystyle\int_{\mathbf{v}\cdot\mathbf{n}<0}(\mathbf{v}\cdot\mathbf{n})\mathcal{M}_{1,\mathbf{u}_{\mathrm w},T_{\mathrm w}}\,\mathrm{d}\mathbf{v}},
\]
with both integrals taken over the corresponding half-ranges of $\Omega_L$.

We use $\mathrm{Kn}=1$ and a uniform $32\times32$ spatial grid and advance to $t=5$. The learned operator is evaluated on a $64^3$ velocity grid in $(-6,6)^3$ through the domain transfer in Section~\ref{sec:velocity-domain-transfer}. Transport uses conservative semi-Lagrangian remapping with CWENO3 reconstruction and a positivity limiter. Each time step consists of a half transport step, a full collision step, and a second half transport step. The collision step uses forward Euler with $\Delta t=0.005$.

The DSMC comparison uses the same spatial grid, $\Delta t=0.003125$, and 2000 particles per cell. We average the moments from 200 independent realizations at each comparison time to obtain the macroscopic quantities.

\begin{figure}[pos=htbp]
\centering
\includegraphics[width=\textwidth]{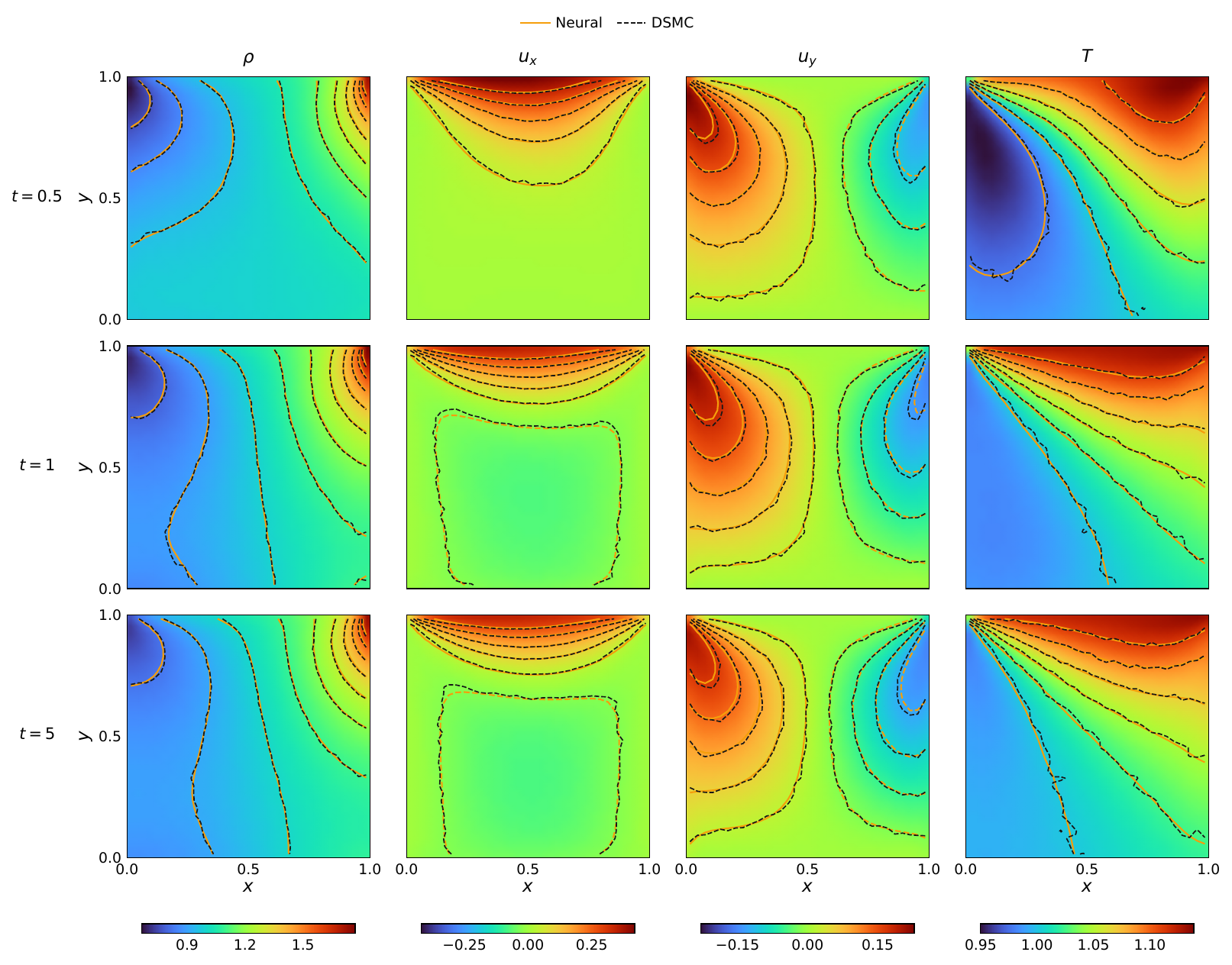}
\caption{Lid-driven cavity at $\mathrm{Kn}=1$. Rows: $t=0.5$, $1$, and $5$. Columns: $\rho$, $u_x$, $u_y$, and $T$. Color fields and orange contours show the solution obtained with the learned operator, and black dashed contours show the DSMC ensemble averages at the same levels.}
\label{fig:cavity-density-temperature}
\end{figure}

Figure~\ref{fig:cavity-density-temperature} compares density, temperature, and velocity at $t=0.5$, $1$, and $5$. The contours obtained with the learned operator closely follow the DSMC ensemble averages. Both computations show density accumulation near the upper-right corner and elevated temperatures beneath the lid, with the warmest region toward the right. The velocity fields show upward motion on the left and downward motion on the right.

\begin{figure}[pos=htbp]
\centering
\includegraphics[width=0.9\textwidth]{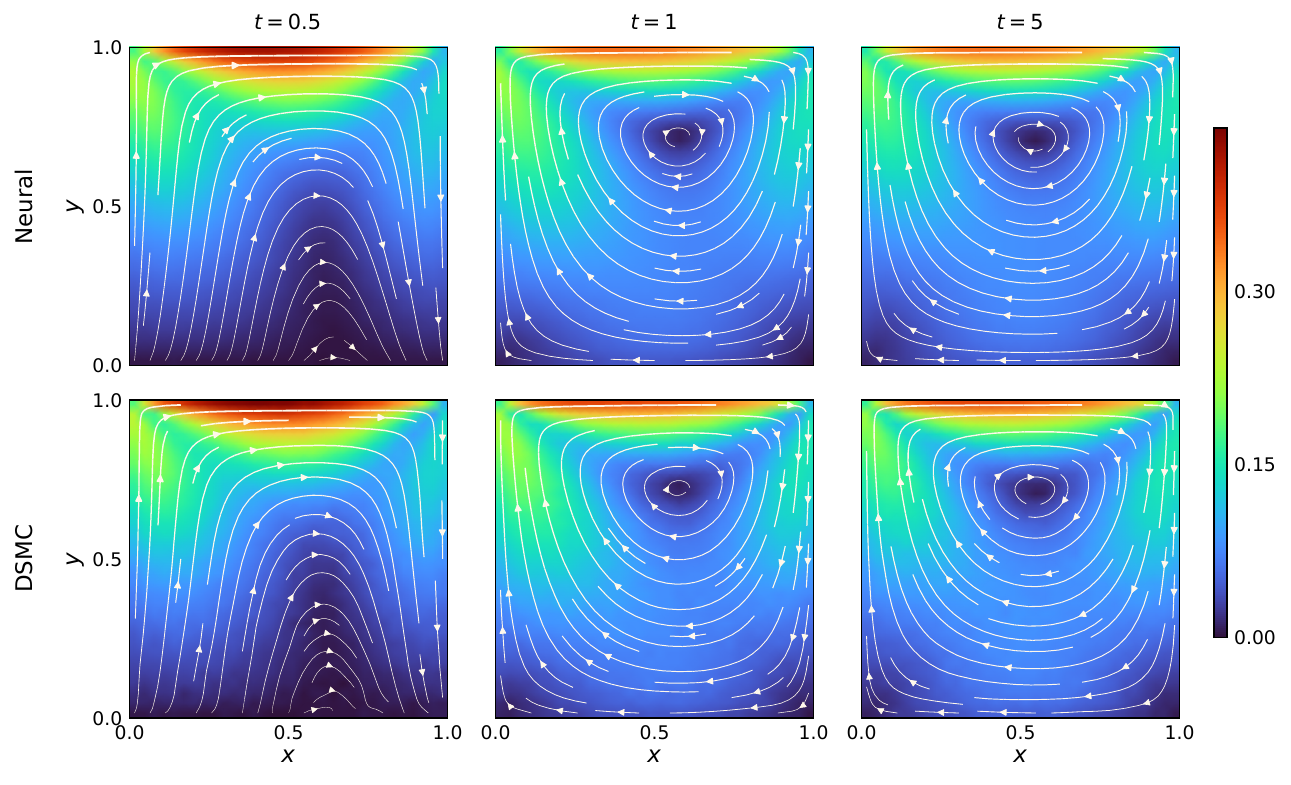}
\caption{Cavity streamlines and speed at $\mathrm{Kn}=1$. Top: learned operator. Bottom: DSMC ensemble average. Columns correspond to $t=0.5$, $1$, and $5$ and share the same color scale for speed.}
\label{fig:cavity-velocity}
\end{figure}

The streamlines in Figure~\ref{fig:cavity-velocity} show how the flow develops. At $t=0.5$, the largest speeds are concentrated beneath the lid. By $t=1$, both computations show a clockwise vortex with a leftward return flow in the lower part of the cavity. The vortex centers remain close to each other at $t=5$.

\begin{figure}[pos=htbp]
\centering
\subfloat[Vertical centerline]{\includegraphics[width=0.42\textwidth]{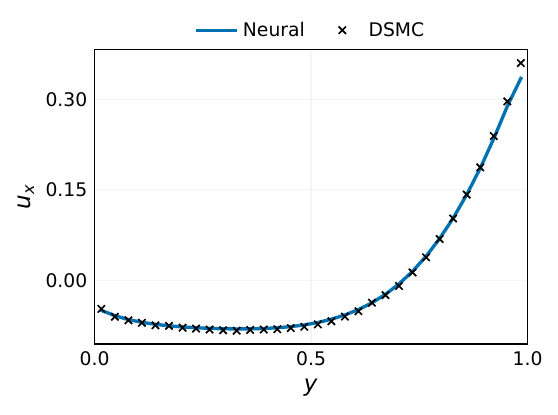}}\quad
\subfloat[Horizontal centerline]{\includegraphics[width=0.42\textwidth]{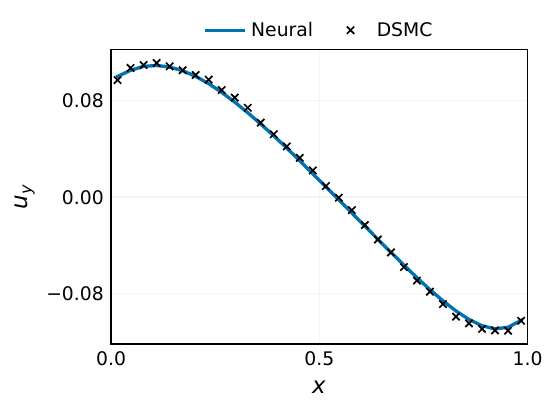}}
\caption{Centerline velocities $u_x(\frac{1}{2},y)$ and $u_y(x,\frac{1}{2})$ at $\mathrm{Kn}=1$ and $t=5$. Solid lines show the solution obtained with the learned operator, and crosses show the DSMC ensemble averages.}
\label{fig:cavity-centerlines}
\end{figure}

At $t=5$, the centerline profiles in Figure~\ref{fig:cavity-centerlines} agree closely with the DSMC ensemble averages. The solution obtained with the learned operator captures the extrema and sign changes of both velocity components, with a small difference in $u_x$ near the moving lid.

The agreement across the time snapshots and centerline profiles shows that repeated collision evaluations with the learned operator capture the development of the cavity flow in a single deterministic simulation.

\section{Conclusion}
\label{sec:conclusion}

We developed a learned Boltzmann collision operator that enforces conservation, entropy dissipation, equilibrium preservation, and quadratic density scaling by construction. Its gain--loss form yields positivity-preserving time discretizations. A positive entropy gap gives exponential relaxation of the homogeneous collision dynamics to the moment-matched equilibrium under the stated assumptions.

The same operator, trained on smooth distributions, captures relaxation from discontinuous data and can be reused in steady Couette flow, normal-shock, and lid-driven cavity solvers. The shock and cavity profiles agree with DSMC ensemble averages.

For the two-stream input, the FP32 learned operator is faster and more accurate than the $p=7$ and $p=9$ spectral evaluations. The accuracy advantage over $p=9$ extends to the evolved distributions in the three-stream and discontinuous-data tests. The learned operator also improves on BGK in collision accuracy, relaxation, and steady flow, and supplies warm starts for spectral refinement.

The experiments use a constant collision kernel. With $C_\omega$ prescribed, an extension to VHS kernels could condition the latent channels on the viscosity index $\omega$ through FiLM~\cite{perez2018film}. Future work will study how collision approximation errors interact with transport discretization and time integration.

\appendix

\section[Generation of the dataset]{Generation of the dataset}
\label{app:dataset-generation}

We sample distributions independently on a uniform $64^3$ grid in $[-10,10)^3$. Each sample starts from a mixture of Gaussian components,
\begin{equation}
  F(\mathbf{v})=\sum_{i=1}^{m}a_i\frac{\exp\left[-\tfrac12(\mathbf{v}-\mathbf{u}_i)^T\Sigma_i^{-1}(\mathbf{v}-\mathbf{u}_i)\right]}{(2\pi)^{3/2}\sqrt{\det\Sigma_i}}.
  \label{eq:dataset-gaussian-mixture}
\end{equation}
We choose $m$ uniformly from $\{1,2,3,4\}$, with positive weights summing to one. For $m>1$, the weights follow a Dirichlet distribution with parameters $0.8$. The centers and covariance orientations are randomized, and the covariance eigenvalues are sampled uniformly on a logarithmic scale over $[0.30,1.80]$. Isotropic covariances recover Maxwellian components, while the sampled covariances allow anisotropy.

With probability $1/2$, we multiply $F$ by
\begin{equation}
  P_\varepsilon(\mathbf{v})=\varepsilon+(1-\varepsilon)\frac{|\mathbf{v}-\mathbf{u}_F|^2}{3T_F},
  \label{eq:dataset-polynomial-modulation}
\end{equation}
where $\mathbf{u}_F$ and $T_F$ are the discrete mean velocity and temperature of $F$. The parameter $\varepsilon$ is zero with probability $0.1$, uniform on $[0,0.1]$ with probability $0.2$, and uniform on $[0.1,1]$ otherwise. This nonnegative modulation introduces depletion near the mean velocity. We denote the resulting distribution, normalized to discrete mass 10, by $g$.

To vary the departure from equilibrium, we form
\begin{equation}
  \widetilde f=(1-\alpha)M_h[g]+\alpha g,
  \label{eq:dataset-equilibrium-blend}
\end{equation}
where $M_h[g]$ is the discrete moment-matched equilibrium from~\eqref{eq:discrete-equilibrium}. This is a BGK-type interpolation that preserves the moments of $g$. Equal numbers of samples are assigned to the strata $[0.05,0.20]$, $[0.20,0.45]$, $[0.45,0.75]$, $[0.75,1]$, and $\{1\}$, with $\alpha$ uniform within each interval. A floor of $10^{-30}$ and normalization to mass 10 are applied when forming both $g$ and the final input $f$.

For each input, we compute the label $\mathcal{Q}_{\mathrm{ref}}[f,f]$ directly by the Fourier spectral method with $\mathcal{B}=1/(4\pi)$, $R=9$, 64 radial quadrature points, and $p=23$. The 10,000 samples are divided into 8,000 training, 1,000 validation, and 1,000 test samples.

\HW{\printcredits}
\section*{\HW{Declaration of competing interest}}
\HW{The authors declare that they have no competing interests.}
\section*{\HW{Funding and acknowledgments}}
\HW{The work of Zhenning Cai was partially supported by the Academic Research Fund of the Ministry of Education of Singapore under grant A-8002392-00-00.}
\section*{Code availability}
The implementation of the learned collision operator and the scripts used to generate the datasets described in Appendix~\ref{app:dataset-generation} are available at \url{https://github.com/shanxue-w/NeuralBoltzmann}.
\section*{\HW{Declaration of generative AI and AI-assisted technologies in the manuscript preparation process}}
\HW{The authors used OpenAI Codex to assist with language editing and \LaTeX{} preparation. The authors reviewed and edited the resulting text and take full responsibility for the content of the manuscript.}

\FloatBarrier
\bibliographystyle{cas-model2-names}
\bibliography{references}

\end{document}